%% file: 0620.tex
\documentstyle[epsf]{amsart}
\theoremstyle{plain}
\newtheorem{Thm}{Theorem}[section]

\newtheorem{Prop}[Thm]{Proposition}

\newtheorem{Lemma}[Thm]{Lemma}
\newtheorem{Conjecture}[Thm]{Conjecture}

\theoremstyle{definition}

\newtheorem{Def*}{Definition}

\newtheorem{Example}[Thm]{Example}
\newtheorem{Cor}[Thm]{Corollary}

\newcommand{\la}{\alpha}
\newcommand{\fkz}{F_{k,Z}(m,n)}
\newcommand{\al}{\alpha}
\newcommand{\be}{\beta}
\newcommand{\ga}{\gamma}

\newenvironment {maplelatex}{}{}
\newenvironment {maplettyout}{}{}
\begingroup\makeatletter\ifx\SetFigFont\undefined%
\gdef\SetFigFont#1#2#3#4#5{%
  \reset@font\fontsize{#1}{#2pt}%
  \fontfamily{#3}\fontseries{#4}\fontshape{#5}%
  \selectfont}%
\fi\endgroup%

\begin{document}

\title[The random assignment problem]{A generalization of\\
the random assignment problem}
\author{Svante Linusson}
\address{Dept of Mathematics\\ Stockholm University\\ 106 91
Stockholm\\ Sweden} \email{linusson@@matematik.su.se }
\author{Johan W\"astlund}
\address{Gustavav\"agen 105\\ S-17831 Eker\"o\\ Sweden}
\email{johan.waestlund@@telia.com}

\date{\today}

\begin{abstract}
We give a conjecture for the expected value of the optimal
$k$-assignment in an $m \times n$-matrix, where the entries are
all exp(1)-distributed random variables or zeros. We prove this
conjecture in the case there is a zero-cost $k-1$-assignment.
Assuming our conjecture, we determine some limits, as $k=m=n\to
\infty$, of the expected cost of an optimal $n$-assignment in an
$n$ by $n$ matrix with zeros in some region. If we take the region
outside a quarter-circle inscribed in the square matrix, this
limit is thus conjectured to be $\pi^2/24$. We give a
computer-generated verification of a conjecture of Parisi for $k=m=n=7$ 
and of a conjecture of Coppersmith and Sorkin for $k\leq 5$. 
We have used the same computer program to verify this conjecture also 
for $k=6$.
\end{abstract}
\maketitle

\section{Introduction}\label{S:Intro}
Suppose we are given an $m$ by $n$ array of nonnegative real
numbers. A \emph{$k$-assignment} is a set of $k$ entries, no two
of which are in the same row or the same column (we ``assign''
rows to columns, or vice versa). The \emph{cost} of the assignment
is the sum of the entries. A $k$-assignment is \emph{optimal} if
its cost is minimal among all $k$-assignments.

We let $F_k(m,n)$ denote the expected cost of the optimal
$k$-assignment in an $m$ by $n$ array of independent exponentially
distributed random variables with mean 1.

Aldous \cite{A} has shown the existence of $\lim_{n\to
\infty}F_n(n,n)=c$. According to a conjecture of M\'ezard and
Parisi \cite{MP85}, $c=\pi^2/6$. Moreover, Parisi \cite{P98} has
conjectured that $F_n(n,n)=1+1/4+1/9+\dots+1/n^2$. This conjecture
has been verified for $n\leq 7$. For $n \le 6$ this can be found
in \cite{svempa}, and the $n=7$ case follows from the calculations
in the appendix to this paper.

In \cite{CS}, a more precise conjecture was formulated.

\begin{Conjecture} \label{C:CS}
$$F_k(m,n)=\sum_{\begin{Sb} i,j\geq 0\\ i+j<k \end{Sb}
}{\frac{1}{(m-i)(n-j)}}.$$
\end{Conjecture}

It is also shown in \cite{CS} that this conjecture is in
accordance with the conjecture of Parisi if $k=m=n$. We have
verified Conjecture \ref{C:CS} for $k\leq 6$, (and arbitrary $m$
and $n$). The cases $k=3,4$ are treated in \cite{svempa} and $k=5$
is given in the appendix. For $k=6$, about 10,000 cases were
needed.

\subsection{The Main Conjecture}
As is shown in \cite{CS}, \cite{svempa}, and in Section
\ref{S:computation}, the calculation of $F_k(m,n)$ can in some
(but not all!) cases be done recursively by reducing the problem
to several assignment problems for matrices which have zeros in
certain positions.

Let $Z$ be a finite set of matrix-positions, that is,
$$Z=\{
(i_1,j_1), \dots, (i_r, j_r)\},$$ for some positive integers
$i_1, j_1,\dots,i_r, j_r$.
Suppose $m\geq \max(i_1,\dots, i_r)$ and
$n\geq \max(j_1,\dots, j_r)$.
We define $F_{k,Z}(m,n)$ to be the expected cost of the optimal
$k$-assignment in an $m$ by $n$ matrix with zeros in the
positions belonging to $Z$, and independent exponential random
variables with mean 1 in all other positions.

We will consider sets of rows and columns in the $m\times
n$-matrix. A set $\la$ of rows and columns is said to {\em cover}
$Z$ if every entry in $Z$ is on either a row or on a column in
$\la$. A covering with $s$ rows and columns will be called an
$s$-covering. $k-1$-coverings will be of particular importance. By
a {\em partial $k-1$-covering} of $Z$, we mean a set of rows and
columns which is a subset of a $k-1$-covering of $Z$. Let
$p_{i,j}$ be the probability that a set of $i$ rows and $j$
columns, chosen from uniform distribution on all such sets, is a
partial $k-1$-covering.

We now present a conjecture for the values of $F_{k,Z}(m,n)$.

\begin{Conjecture}[Main Conjecture] \label{C:main}
\begin{equation} \label{eq:main}
F_{k,Z}(m,n) = \sum_{\begin{Sb} i,j\geq 0\\ i+j<k
\end{Sb}}{\frac{p_{i,j}}{(m-i)(n-j)}}
\end{equation}
\end{Conjecture}

\begin{Example}
When $Z$ is empty, Conjecture \ref{C:main} clearly specializes to
Conjecture \ref{C:CS}.
\end{Example}

\begin{Example}\label{E:onezero}
When $k=2$ and $Z$ has just one element, then $p_{1,0}=1/m$,
$p_{0,1}=1/n$ and $p_{0,0}=1$, so according to Conjecture
\ref{C:main} we get
$$F_{2,Z}(m,n)=\frac{1}{mn}+\frac{1/m}{(m-1)n}+
\frac{1/n}{m(n-1)}=\frac{-1}{mn}+\frac{1}{(m-1)n}+\frac{1}{m(n-1)},$$
which is also the correct value, see Example \ref{E:k=2} below.
\end{Example}

In fact, we prove in Section \ref{S:proof} that Conjecture
\ref{C:main} holds under the assumption that $Z$ contains a
$k-1$-assignment, which is our main theorem.

Let $P$ be the poset of all intersections of $k-1$-coverings of $Z$,
ordered by reverse inclusion, and with an artificial minimal element
$\hat{0}$ corresponding to the empty intersection. We prove in Section
\ref{S:proof} that the following is an equivalent form of the
Conjecture \ref{C:main}:

\begin{equation} \label{eq:reformulation}
F_{k,Z}(m,n) = \sum_{\la\in P\backslash \hat 0}
{\frac{-\mu_P(\al)}{(m-i_\al)(n-j_\al)}},
\end{equation}
where $\mu_P(\alpha)$ is the M\"obius function of the interval
$(\hat{0},\al)$ in $P$ and $i_\al,j_\al$ are the number of rows and
columns in $\al$ respectively.

\smallskip
In Section \ref{S:computation} we describe our tools to
recursively compute $\fkz$, which we use both for the theoretic
proofs and for the algorithmic computational results in the
Appendix. In Section \ref{S:consequences} we discuss some
consequences of Conjecture \ref{C:main}. One such consequence is
that we always can write
$$\fkz=\sum_{i+j<k}{\frac{b_{i,j}}{(m-i)(n-j)}},$$ for some
integers $b_{i,j}$ independent of $m$ and $n$. We also prove that
Conjecture \ref{C:main} implies a conjecture by Olin \cite {Olin}
on the probability that the smallest element in a row is used in
the optimal assignment. In Section \ref{S:asymptotics} we prove
that some asymptotic results follows from Conjecture \ref{C:main}.

\noindent {\bf Acknowledgement:} We thank Mireille
Bousquet-M\'elou and Gilles Schaeffer for introducing us to the
problem.

\section{Computation of $F_{k,Z}(m,n)$}\label{S:computation}
If $A$ is a (random) $m$ by $n$ matrix, and $u$ is an assignment,
then we let $\Psi_u(A)$ be the cost of $u$, that is, $$
\Psi_u(A)=\sum_{(i,j)\in u}{A_{i,j}}.$$ If $k$ is an integer not
greater than $m$ or $n$, we let $\Psi_k(A)$ be the cost of the
optimal $k$-assignment in $A$.

We say that a set $\la$ of rows and columns is an \emph{optimal
covering} of a set $Z$, if $\la$ covers $Z$, and $\la$ has minimal
cardinality among all coverings of $Z$.
Our calculations of $F_{k,Z}(m,n)$ are based on the following
theorem.

\begin{Thm}\label{T:recursion}
Let $k\le m,n$ be positive integers.
Let $A$ be an $m$ by $n$ matrix, with zeros in a certain set $Z$
of positions, and (possibly random) positive values outside $Z$.
Let $\al$
be an optimal covering of $Z$. Let $\epsilon$ be a positive real
(possibly random) number, not greater than any of the entries not
covered by $\al$. Let $A'$ be the matrix obtained from $A$ by
subtracting $\epsilon$ from every entry not covered by $\al$, and
adding $\epsilon$ to every doubly covered entry of $A$. Then
$$\Psi_k(A)=\epsilon\cdot(k-\left| \al\right|) + \Psi_k(A').$$
\end{Thm}

To prove this theorem, we need the following lemma, which we prove
in Section \ref{S:lemmas}.

\begin{Lemma}\label{L:subset}
If $\al$ is an optimal covering of a subset of $Z$, then there is
an optimal $k$-assignment which intersects every row and every
column of $\al$.
\end{Lemma}

\begin{pf*}{Proof of Theorem \ref{T:recursion}}
Let $u$ be an optimal $k$-assignment of $A'$ intersecting every
row and column of $\al$. Such an optimal assignment exists by
Lemma \ref{L:subset}. Note that subtracting $\epsilon$ from every
entry not covered by $\al$, and adding $\epsilon$ to every doubly
covered entry is the same thing as first subtracting $\epsilon$
from all entries not covered by the rows of $\al$, and then adding
$\epsilon$ to all the entries covered by the columns of $\al$. Let
$i$ and $j$ be the number of rows and columns in $\al$,
respectively. Then
\begin{multline}\label{eq:psi}
\Psi_u(A)=\Psi_u(A')+\epsilon\cdot(k-i)-\epsilon\cdot j=\\
=\epsilon\cdot(k-\left|
\al\right|)+\Psi_u(A')=\epsilon\cdot(k-\left|
\al\right|)+\Psi_k(A').
\end{multline}
For every $k$-assignment $v$ we have $$ \Psi_v(A)\geq
\epsilon\cdot(k-\left| \al\right|)+\Psi_v(A')\geq
\epsilon\cdot(k-\left| \al\right|)+\Psi_k(A'). $$ Hence $u$ is an
optimal $k$-assignment also in $A$, and $\Psi_k(A)$ is given by
\eqref{eq:psi}.
\end{pf*}

In the typical use of the theorem, $\epsilon$ is the minimum of
all the non-covered entries, which gives a new zero when we
subtract $\epsilon$. The special case of this theorem where there
are no doubly covered entries or all such entries are known not to
be in the optimal $k$-assignment, was proved and used in \cite{CS}
and \cite{svempa} to prove Conjecture \ref{C:CS} for $k\le 4$.

Before proving Lemma \ref{L:subset}, we give an example of how
Theorem \ref{T:recursion} is used. First we need to state some
well-known properties of the exponential distribution. We say that
a random variable is $exp(l)$-distributed if it is exponentially
distributed with mean $l$.

\begin{Lemma}
Let $X_1,\dots,X_t$ be independent exponential random variables,
with mean $a_1, a_2,\dots, a_t$ respectively. Let
$Y=\min\{X_1,\dots,X_t\}$. Then
\begin{enumerate}
\item{} $Y$ is exponentially distributed with mean $$\frac{1}{1/a_1+\dots+1/a_t}$$
\item{} if $Y=X_i$, then $X_j-Y$ is again
$exp(a_j)$-distributed, for $i\neq j$, and $\{X_j-Y: j\neq i\}$
are independent, and independent of $Y$.
\item{} if $I$ is the random variable defined by $Y=X_I$, then $I$ is
uniquely defined with probability 1 and
$$Pr(I=i)=\frac{1/a_i}{1/a_1+\dots+1/a_t}.$$
\end{enumerate}
\end{Lemma}

\begin{Example}\label{E:k=2}
Let us show how this is used to compute $F_2(m,n)$. First we take
$\epsilon$ to be the minimum of all elements of $A$, which by
symmetry we can assume to be $a_{1,1}$. We now subtract $a_{1,1}$
from all entries in $A$ to obtain the matrix $A'$ with a zero in
position $(1,1)$ and exp(1)-distributed random variables
elsewhere. The covering $\al$ is here empty and $\epsilon$ is
exp($1/mn$)-distributed. From Theorem \ref{T:recursion} we get
\[
F_2(m,n)=\frac{2}{mn}+F_{2,Z'}(m,n),
\] where $Z'=\{(1,1)\}$.

To compute $F_{2,Z'}(m,n)$, we have to find an optimal covering of
$Z'$ and we can choose $\al'= \{r_1\}$, the first row, see Figure
\ref{fg:k=2}. We again take the minimum of the entries that are not
covered by $\al'$, which means that $\epsilon'$ is $\exp(1/(m-1)n)$.
{}From Theorem \ref{T:recursion} we get
\[
F_{2,Z'}(m,n)=\frac{1}{(m-1)n}+\frac{1}{n}F_{2,Z''}(m,n)+
\frac{n-1}{n}F_{2,Z'''}(m,n),
\]
where $Z''=\{(1,1),(2,1)\}$ and $Z'''=\{(1,1),(2,2)\}$. The
minimum is equally likely to occur in any of the entries that is
not in the first row. Here $Z''$ represents the case that the
minimum occurs somewhere in the first column, while $Z'''$
represents the case that the minimum occurs in another column.
Hence we get only two different non-isomorphic cases.

$F_{2,Z'''}(m,n)=0$ because $Z'''$ contains a 2-assignment of
zeros. Finally we use Theorem \ref{T:recursion} on $Z''$. The only
optimal covering is the first column, $\al''=\{c_1\}$. The minimum
of the other entries is $\exp(1/m(n-1))$, which we use as
$\epsilon''$. Note that a new zero in another column than the
first will produce a 2-assignment together with one of the zeros
in $Z''$. We get $$F_{2,Z''}(m,n)=\frac{1}{m(n-1)}.$$ Substituting
back we get

$$F_{2,Z'}(m,n)=\frac{-1}{mn}+\frac{1}{(m-1)n}+\frac{1}{m(n-1)},$$
as in Example \ref{E:onezero}, and

$$F_{2,Z}(m,n)=\frac{1}{mn}+\frac{1}{(m-1)n}+\frac{1}{m(n-1)},$$
in accordance with Conjecture \ref{C:main}.

\begin{figure}[htb]
\centerline{\epsfysize=34mm\epsfbox{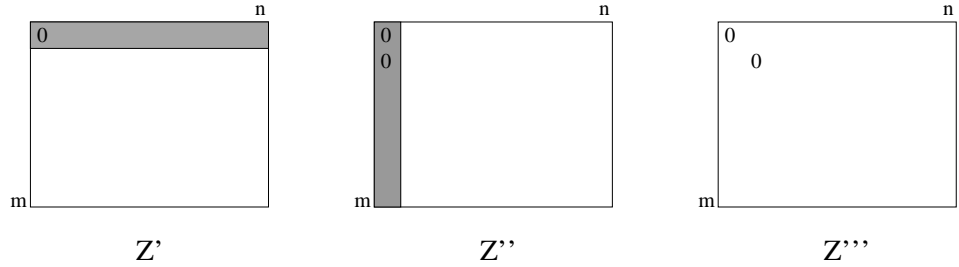}}
\caption{\label{fg:k=2} The cases needed when computing $F_2(m,n)$.
The shaded row and column are the coverings. }
\end{figure}
\end{Example}

\begin{Example}\label{E:firstbad}
Let us also give an example where there is a doubly covered entry.
Let $k=5$ and $Z=\{(1,2),(1,3),(2,1),(3,1)\}$. This is the
simplest case where the recursion will result in a matrix with an
entry which is neither $\exp(1)$ nor zero.

The only optimal covering $\al=\{r_1,c_1\}$ consists of row one
and column one. As usual we let $\epsilon$ be the minimum of the
non-covered entries, which here is $\exp(1/(m-1)(n-1))$. We get

\[
F_{5,Z}(m,n)=\frac{3}{(m-1)(n-1)}+ \frac{1}{(m-1)(n-1)}
\sum_{(i,j)\in[2,m]\times [2,n]}\Psi_5(A_{i,j}'),
\]
where $A_{i,j}'$ has zeros in $Z\cup \{(i,j)\}$ and an entry which
is the sum of two r.v.'s, $a_{1,1}=\exp(1)+\exp(1/(m-1)(n-1))$.
There are actually four different cases here, and as is shown in
the appendix it is possible to continue the recursion in each case
to compute $F_{5,Z}(m,n)$ and finally $F_{5}(m,n)$.

\begin{figure}[htb]
\centerline{\epsfysize=34mm\epsfbox{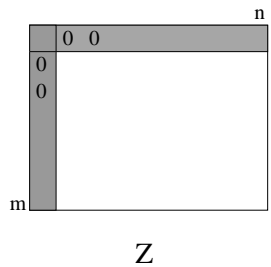}}
\caption{\label{fg:firstbad}
The smallest ``problem case''. The shaded row and
column is the covering.
}
\end{figure}
\end{Example}

Note that in Section \ref{S:prob} we show a way to avoid getting sums
of random variables as entries in this particular case. We can thus
reduce the number of ``problem cases'', that is cases when the
recursion in Theorem \ref{T:recursion} gives us entries different from
$\exp(1)$ and zero. But we have not found any way to avoid them 
completely.

\subsection{Proof of Lemma \ref{L:subset}} \label{S:lemmas}

We begin by citing a famous theorem of D. K\"onig:

\begin{Thm}
If a set of matrix positions cannot be covered by $k-1$ rows and
columns, then it contains a $k$-assignment.\qed
\end{Thm}

\begin{Lemma}\label{L:rowsused}
Let $k\le m,n$ be positive integers.
Let $A$ be an $m$ by $n$ matrix, with zeros in a certain set $Z$
of positions, and (possibly random) positive values outside $Z$.
Suppose that $Z$ does not contain a $k+1$-assignment. Let $\al$ be an
optimal covering of $Z$. Then every row and every column of $\al$
contains an element of every optimal $k$-assignment of $A$.
\end{Lemma}

\begin{pf}
Suppose that $\al$ contains rows $1,\dots, r$ and no other rows.
Since $\al$ is optimal it follows by K\"onig's theorem that there
is an $r$-assignment $v$ in $Z$ containing no element from the
columns in $\al$. Suppose (for a contradiction) that there is an
optimal $k$-assignment $u$ which does not use row 1.

We now construct a sequence of matrix positions (all with zeros)
as follows: Let
$v_0$ be the element of $v$ which is in the first row. Suppose
that we have defined $v_0,\dots,v_h$ and $u_1,\dots,u_h$. Then if
there is an element of $u\cap Z$ in the same column as $v_h$, let
this element be $u_{h+1}$, and let $v_{h+1}$ be the element of $v$
which is in the same row as $u_{h+1}$. Since $u$ does not
contain any element from the first row, the sequences $v_0, v_1,
v_2,\dots$ and $u_1, u_2,\dots$ cannot contain any repetitions of
the same element. Hence the sequence must end with an element
$v_h$ such that there is no element of $u\cap Z$ in the same
column.

We now consider two cases. Suppose first that
no element in the column of $v_h$ belongs to $u$. Since $Z$ does
not contain a $k+1$-assignment, the cardinality of $\al$ is at most
$k$. Each row and column of $\al$ covers at most one element of $u$,
and row 1 does not cover any element of $u$. Consequently there is
an element $(i,j)$ of $u$ which is not covered by $\al$, hence does
not belong to $Z$. Then $u\setminus\{u_1,\dots,
u_{h-1},(i,j)\}\cup\{v_0,\dots,v_{h}\}$ is a $k$-assignment of
smaller cost than $u$, a contradiction. If on the other hand there
is an element $(i',j')$ of $u$ in the first column, then
$u\setminus\{u_1,\dots, u_{h-1},(i',j')\}\cup\{v_0,\dots,v_{h}\}$
is a $k$-assignment of smaller cost than $u$. This contradiction
proves the lemma.
\end{pf}

\begin{pf*}{Proof of Lemma \ref{L:subset}}
This follows by a continuity argument. Let the values in the
positions of $Z$ not covered by $\al$ be $\delta$, and let
$\delta$ tend to zero. Since there are only finitely many
$k$-assignments, there has to be an assignment which is optimal
for all sufficiently small $\delta$, which by Lemma
\ref{L:rowsused} has to intersect every row and column of $\al$.
This assignment is also optimal for $\delta =0$.
\end{pf*}

Note that it is not true that every optimal $k$-assignment has to
intersect the rows and columns of an optimal covering of a subset
of $Z$.

\subsection{Superfluous matrix elements}
It can be of great help to know that a certain entry of the matrix
cannot occur in any optimal assignment. Sufficient conditions to draw
this conclusion were given in Lemma 26 of \cite{CS} and Lemma 9 of
\cite{svempa}. The following lemma gives a necessary and sufficient
condition. The condition here is in fact equivalent to the condition
given in \cite{svempa}.

\begin{Thm} \label{T:insertzero}
Let $k\le m,n$ be positive integers. Let $A=(a_{i,j})$ be an $m$
by $n$ matrix. Suppose that a set $Z$ of entries are zero, and the
remaining entries are random variables that can take arbitrary
positive values. Suppose that $(i,j)\notin Z$. If every set of
$k-1$ rows and columns covering $Z$ also covers $(i,j)$, then
$(i,j)$ cannot belong to an optimal $k$-assignment. Conversely, if
there is a $k-1$-covering of $Z$ which does not cover $(i,j)$,
then it is possible to assign positive values to the matrix
entries not in $Z$ in such a way that $(i,j)$ belongs to a unique
optimal $k$-assignment.

\end{Thm}

\begin{pf}
By the theorem of K\"onig, a set contains a $k$-assignment if and
only if it cannot be covered by fewer than $k$ rows and columns.
The minimal cost of a $k$-assignment is of course the same as the
minimal cost of a set containing a $k$-assignment. In trying to
minimize the cost of such a set, it suffices to consider sets
containing $Z$, since the entries in $Z$ do not increase the
cost.

Suppose that $(i,j)$ is covered whenever
$Z$ is covered with fewer than $k$ rows and columns. Then
from any set containing $Z$, which cannot be covered with
fewer than $k$ rows and columns, we can delete $(i,j)$ and obtain
another such set of smaller cost. This shows that $(i,j)$ can
never belong to an optimal $k$-assignment.

Suppose on the other hand that we can cover $Z$ with a set $\al$
of $k-1$ rows and columns, without covering $(i,j)$. Then we
assign the value 1 to $a_{i,j}$, and to the entries not in $Z$
which are covered by $\al$, and values larger than $k$ to all
other entries not in $Z$. Since $m$ and $n$ are not smaller than
$k$, the set of entries which are covered by $\al$ cannot be
covered in any other way by $k-1$ rows and columns. Hence the set
of entries which are $\leq 1$ cannot be covered with $k-1$ rows
and columns. By K\"onig's theorem there is a $k$-assignment of
cost at most $k$. However, since every $k$-assignment must use at
least one entry not covered by $\al$, every $k$-assignment of cost
at most $k$ must use $a_{i,j}$. Hence $a_{i,j}$ belongs to the
unique optimal $k$-assignment.
\end{pf}

In this situation, the fact that $F_{k,Z}(m,n)=F_{k,Z\cup
\{(i,j)\}}(m,n)$ can sometimes simplify computation. It is
possible that several elements have the property of $(i,j)$ in the
above theorem. Indeed, as the following theorem shows, it is
sometimes possible to dispose of an entire row or column. If $r$ is a
row, we let $Z\backslash r$ denote the set where all positions of $Z$
in row $r$ has been removed and similarly for a column.

\begin{Thm} \label{T:deleterow}
Suppose that a row $r$ is used in every covering of $Z$ by $k-1$
rows and columns (in particular, this is the case if at least $k$
elements in $r$ belong to $Z$). Let $Z'=Z\setminus r$. Then
\begin{equation}\label{eq:deleterow}
F_{k,Z}(m,n)=F_{k-1,Z'}(m-1,n).
\end{equation}
Similarly, if a column $c$ belongs to every covering of $Z$ with
fewer than $k$ rows and columns, and $Z''=Z\setminus c$, then
$F_{k,Z}(m,n)=F_{k,Z''}(m,n-1)$.
\end{Thm}

\begin{pf}
If $r$ is used in every covering of $Z$ by $k-1$ rows and columns,
then by Theorem \ref{T:insertzero}, no non-zero element in row
$r$, no matter how small, is used in any optimal $k$-assignment.
Hence we do not change the cost of the optimal $k$-assignment if
we replace all non-zero entries in row $r$ by zeros. Let $A'$ be
the matrix obtained from $A$ by deleting the row $r$. Every
$k$-assignment of $A$ must of course contain (at least) a
$k-1$-assignment of $A'$. It is not a problem if $Z$ contains an entry
in row $m$, since we can always permute the rows in $Z$ so row $m$
becomes empty instead of row $r$.

Conversely, if row $r$ consists of only
zeros, then every $k-1$-assignment in $A'$ can be extended to a
$k$-assignment in $A$ of the same cost. Hence
$\Psi_k(A)=\Psi_{k-1}(A')$. Passing to expected values, we obtain
\eqref{eq:deleterow}.

The corresponding statement for columns is of course equivalent.
\end{pf}

\subsection{The probability that a certain element belongs to the
optimal assignment}\label{S:prob}

In this section we prove a theorem relating the expected values of
certain assignments to the probability that a certain element
belongs to such an assignment. We show how this implies
a second recursion which can be used to avoid
some ``problem cases'' in the computation of $\fkz$.

\begin{Thm} \label{T:e-d}
Let $k\le m,n$ be positive integers. Let $A=(a_{i,j})$ be an
$m\times n$-matrix, where a set $Z$ of the entries are zero, and
the others are independent exponentially distributed random
variables with mean 1. Suppose that $(i,j)\notin Z$, and let
$Z'=Z\cup \{(i,j)\}$. Then the probability that $(i,j)$ belongs to
an optimal $k$-assignment in $A$ is $F_{k,Z}(m,n)-F_{k,Z'}(m,n)$.
\end{Thm}

\begin{pf}
We condition on the values of all matrix elements except
$a_{i,j}$. Let $X$ be the value of the minimal $k$-assignment in
$A$ which does not use $a_{i,j}$. Let $Y$ be the value of the
minimal $k-1$-assignment in $A$ which does not use row $i$ or
column $j$. The optimal $k$-assignment in $A$ either contains
$(i,j)$ and has value $Y+a_{i,j}$, or does not contain $(i,j)$ and
has value $X$. Hence the probability that $a_{i,j}$ belongs to the
optimal $k$-assignment in $A$ is equal to the probability that
$a_{i,j}<X-Y$.

We wish to show that this is equal to
$$F_{k,Z}(m,n)-F_{k,Z'}(m,n)=\text{E}(\max(0,\min(X-Y, a_{i,j})).$$

If $X\leq Y$, then both are zero. If $X>Y$, then we let
$\delta=X-Y$. Then the probability that $a_{i,j}$ is used in the
optimal $k$-assignment in $A$ is $1-e^{-\delta}$. We compare this
to $$\text{E}(\min(\delta, a_{i,j}))=\delta
e^{-\delta}+\int_{0}^{\delta}{te^{-t}\,dt} = \delta
e^{-\delta}+1-(\delta+1)e^{-\delta}=1-e^{-\delta}.$$ This proves
the theorem.
\end{pf}

{}From this we can get another generalization of the basic
recursion used in \cite{CS} and \cite{svempa}, which is of some
help when all the doubly covered entries are nonzero. However,
none of the results in this paper depends on this recursion.

\begin{Cor}\label{C:recursion2}
Let $k\le m,n$ be integers. Let $Z\subset [1,m]\times [1,n]$, be a
set of zeros. Let $\al$ be an optimal covering of $Z$ and let $R$
be the non-covered entries. Also let the random variable $S$ be
the number of positions not covered by $\al$ that are in the
optimal assignment. Then
\[F_{k,Z}(m,n)=\frac{E[S]}{|R|}+
\frac{1}{|R|}\sum_{(i,j)\in R}F_{k,Z\cup(i,j)}(m,n).
\]
\end{Cor}
\begin{pf}
By linearity of the expected value we have $$E[S]=\sum_{(i,j)\in
R} (F_{k,Z}(m,n)-F_{k,Z\cup(i,j)}(m,n).$$
\end{pf}

\begin{Example}
Let us give an example of how Corollary \ref{C:recursion2} can be used
to simplify the calculations of $\fkz$. As in Example \ref{E:firstbad}
we take $k=5$ and $Z=\{(1,2),(1,3),(2,1),(3,1)\}$. The optimal
covering $\al$ is the first row and the first column. The expected
number of non-covered elements used in the optimal k-assignment is
$3+$ the probability that $a_{1,1}$ is in the optimal assignment which
is $F_{k,Z}(m,n)-F_{k,Z\cup \{(1,1)\}}(m,n)$. We get

\begin{multline}
F_{k,Z}(m,n)=\frac{3+F_{k,Z}(m,n)-F_{k,Z\cup
\{(1,1)\}}(m,n)}{(m-1)(n-1)} +\\
+\frac{1}{(m-1)(n-1)}\sum_{(i,j)\in [2,m]\times[2,n]}
F_{k,Z\cup \{(i,j)\}}(m,n).
\end{multline}

This can then be rewritten as
\[F_{k,Z}(m,n)=\frac{3-F_{k,Z\cup \{(1,1)\}}(m,n)+
\sum_{(i,j)\in [2,m]\times[2,n]}
F_{k,Z\cup \{(i,j)\}}(m,n)}{(m-1)(n-1)-1},
\]
which enables us to stay within the set of matrices with all
entries exp(1) or zero. However, this recursion does not help in
the case $$Z=\{(1,1),(1,2),(1,3),(2,1),(3,1)\}.$$ For this we
would need to know the probability that the zero in position
$(1,1)$ is in the optimal k-assignment. We know of no simple way
to calculate this, and in general it might not even be
well-defined, since in some cases with positive probability there
are several optimal $k$-assignments using the same non-zero
entries, but different zeros.
\end{Example}

\section{Proof of the Main Theorem}\label{S:proof}
In this section we give a proof of Conjecture \ref{C:main} under
the assumption that $Z$ contains a $k-1$-assignment.

Let $k$, $m$, $n$ and $Z$ be as usual. Let $B(m,n)$ be the Boolean
algebra of all sets of rows and columns in $A$, ordered by reverse
inclusion. We define a function $f=f_{k,Z}$ on $B(m,n)$ by

$$ f(\la)=\begin{cases} 1, \,\text{if $\la$ is a partial
$k-1$-covering}\\ 0, \,\text{otherwise}
\end{cases}$$

We define $g=g_{k,Z}$ by $$\sum_{\be\leq \la}{g(\be)}=f(\la).$$

\begin{Lemma} \label{L:inducedposet}
Let $P=P_{k,Z}(m,n)$ be the sub-poset of $B(m,n)$ consisting of
all intersections of $k-1$-coverings, together with an artificial
minimal element. Then $$g(\la)=\begin{cases} -\mu_P(\la),\,
\text{if $\la\in P$}\\ 0,\, \text{if $\la\notin P$}\end{cases}$$
\end{Lemma}

\begin{pf}
Suppose that $\la$ is minimal contradicting the statement. Then
$\la\notin P$ and $g(\la)\neq 0$. Let $\la'$ be the intersection
of all atoms of $P$ which are smaller than $\la$, that is, which
are supersets of $\la$. Then $\la'\in P$, and for every $\be\in
P$, if $\be\leq \la$, then $\be\leq \la'$. Hence
$\sum_{\be<\la}{g(\be)}=0$, and consequently, $g(\la)=0$. This is
a contradiction.
\end{pf}

For $\la \in B(m,n)$, we let $i_\la$ and $j_\la$ be the number of rows
and columns in $\la$, respectively.

As we will prove, a reformulation of our main conjecture is the
following:

\begin{Conjecture} \label{C:mainmoeb}
\begin{equation} \label{eq:gsum}
F_{k,Z}(m,n)=\sum_{\la\in B(m,n)}{\frac{g(\la)}{(m-i_\la)(n-j_\la)}}.
\end{equation}
\end{Conjecture}

\begin{Example}\label{E:poset}
Let $k=3$ and $Z=\{(1,1)\}$. We will compute $F_{3,Z}(m,n)$ using
Conjecture \ref{C:mainmoeb} and Lemma \ref{L:inducedposet}. The
$k-1$-coverings of $Z$ are all sets of two rows and columns
including row one or column one. The poset $P=P_{3,Z}(m,n)$ and
the M\"obius function is shown in Figure \ref{fg:poset}. The sets
with an index $i$, e.g. $r_1r_i$, represents all $m-1$ sets with
$2\le i \le m$ and similarly the sets with an index $j$ represents
all $n-1$ sets with $2\le j\le n$. We get

\begin{multline}
F_{3,Z}(m,n)=(m-1)\frac{1}{(m-2)n}+(n-1)\frac{1}{(m-1)(n-1)}+
\frac{1}{(m-1)(n-1)}+\\
+(m-1)\frac{1}{(m-1)(n-1)}+(n-1)\frac{1}{m(n-2)}+
\frac{-m-n+2}{(m-1)n}+\\+(m-1)\frac{-1}{(m-1)n}+
(n-1)\frac{-1}{m(n-1)}+ \frac{-m-n+2}{m(n-1)}+\frac{m+n-2}{mn}=\\
=\frac{1}{(m-2)n}+\frac{1}{m(n-2)}+\frac{1}{(m-1)(n-1)}
+\frac{1}{m(n-1)}+\frac{1}{(m-1)n}+\frac{-2}{mn},
\end{multline}
which is the correct answer.

\begin{figure}[htb]
\centerline{\epsfysize=60mm\epsfbox{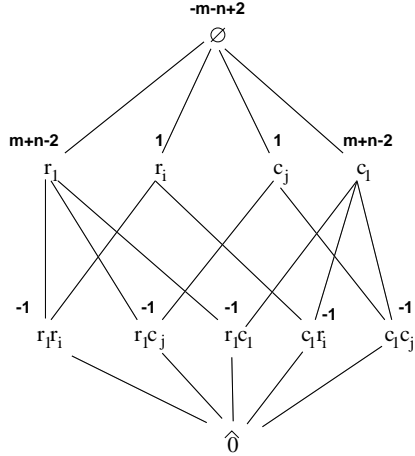}}
\caption{\label{fg:poset}
The poset $P_{3,Z}(m,n)$ in Example \ref{E:poset}. The M\"obius
function is written in boldface.
}
\end{figure}
\end{Example}

\medskip
We now state a lemma, which will also be useful later.
It should be known, but for completeness we give a proof.

\begin{Lemma} \label{L:sum}
If $m$ and $r$ are positive integers, then $$
\sum_{i=0}^{r}{\frac{(-1)^i\binom{r}{i}}{m-i}}=\frac{(-1)^r}{m\binom
{m-1}{r}}.$$
\end{Lemma}

\begin{pf}
We have
\begin{multline} \label{eq:lm}
\sum_{i=0}^{r}{\frac{(-1)^i\binom{r}{i}}{m-i}}=
\sum_{i=0}^{r}{\int_{0}^{1}{(-1)^i\binom{r}{i}x^{m-1-i}\, dx}}=\\
=\int_0^1{\left(
\sum_{i=0}^{r}{(-1)^i\binom{r}{i}x^{m-1-i}}\right)\, dx}=\\
=(-1)^r\int_0^1{x^{m-1-r}
\sum_{i=0}^{r}{(-1)^{r-i}\binom{r}{r-i}x^{r-i}}\, dx}=\\
=(-1)^r\int_0^1{x^{m-1-r}(1-x)^{r} \, dx}.
\end{multline}

This integral can be interpreted as a probability. Let
$Y_1,\dots,Y_{m}$ be independent random variables uniformly
distributed in $[0,1]$. Then the integral in \eqref{eq:lm} represents
the probability that $Y_i<Y_{m-r}$ for $i=1,\dots,m-r-1$, and
$Y_j>Y_{m-r}$ for $j=m-r+1,\dots,m$. If we sort the $m$ numbers
$Y_1,\dots,Y_{m}$ according to size, every permutation occurs with
equal probability $1/m!$. The number of permutations such that $Y_i,
i=1,\dots,m-r-1$ occur before $Y_{m-r}$, and $Y_j$, $j=m-r+1,\dots,m$
all occur after $Y_{m-r}$, is $r!(m-r-1)!$. Therefore, the probability that
$Y_{i}<Y_{m-r}<Y{j}$ for all $i=1,\dots,m-r-1$ and $j=m-r+1,\dots,m$, is
$$ \frac{r!(m-r-1)!}{m!}=\frac{1}{m\binom{m-1}{r}}.$$ This proves the
lemma.
\end{pf}

\begin{Prop}
Conjecture \ref{C:main} $\Leftrightarrow$ Conjecture
\ref{C:mainmoeb}.
\end{Prop}

\begin{pf}
Let $C$ be the set of partial $k-1$-coverings. For $\be\in C$, we
define
\[
f_\be(\al)=
\begin{cases}
1,\, \text{if $\al=\be$}\\ 0,\, \text{otherwise}
\end{cases}
\]

Obviously, $f=\sum_{\be\in C}{f_\be}$. Moreover,
\[
g=\sum_{\be\in C}{g_\be},
\]
where $g_\be$ is defined by
$$\sum_{\ga\leq \la}{g_\be(\ga)}=f_\be(\la).$$

We have $g_\be(\la)=\mu(\be,\la)$, the M\"obius function of the interval
$(\be,\la)$ in the Boolean algebra over the set of rows and columns
(ordered by reverse inclusion). This interval is itself a Boolean
algebra. Hence
\[
g_\be(\la)=
\begin{cases}
(-1)^{|\be|-|\la|},\, \text{if $\la\subseteq \be$}\\ 0,\, \text{otherwise}
\end{cases}
\]

For the right hand side of \eqref{eq:gsum}, we have
\begin{multline} \label{eq:probsum}
\sum_{\la\in B(m,n)}{\frac{g(\la)}{(m-i_\la)(n-j_\la)}}=
\sum_{\be\in C}{\sum_{\la\geq
\be}{\frac{(-1)^{|\be|-|\la|}}{(m-i_\la)(n-j_\la)}}}
=\\
=\sum_{\be\in
C}{\left(\sum_{i=0}^{i_\be}{(-1)^i\frac{\binom{i_\be}{i}}{m-i}}\right)
\left(\sum_{j=0}^{j_\be}{(-1)^j\frac{\binom{j_\be}{j}}{n-j}}\right)}=\\
=[\text{by Lemma
\ref{L:sum}}]=\sum_{\be\in
C}{\frac{1}{m\binom{m-1}{i_\be}n\binom{n-1}{j_\be}}}=\\
=\sum_{\be\in C}{\frac{1}{\binom{m}{i_\be}\binom{n}
{j_\be}(m-i_\be)(n-j_\be)}}.
\end{multline}

The factor
\[
\frac{1}{\binom{m}{i_\be}\binom{n}{j_\be}}\]
can be
interpreted as the probability of $\be$ under the uniform
probability measure on all sets of $i_\be$ rows and $j_\be$ columns.
Hence \eqref{eq:probsum} equals \eqref{eq:main}.
\end{pf}

We need the following lemma.

\begin{Lemma}
Suppose $Z$ contains a $k-1$-assignment $u$. Let $Z'=Z\cup
\{(i,j)\}$, where $(i,j)$ lies in a column not used by $u$. Then
$$ g_{k,Z'}(\al)=
\begin{cases}
g_{k,Z}(\al),\, \text{if $\al$ contains row $i$}\\
0,\, \text{otherwise}
\end{cases}
$$
\end{Lemma}

\begin{pf}
It follows from the assumptions that in order to cover $Z'$ with
$k-1$ rows and columns, one must use row $i$. This implies that
$P_{k,Z'}(m,n)$ consists of precisely those elements of
$P_{k,Z}(m,n)$ which contain row $i$. Since this means that
$P_{k,Z'}(m,n)$ is an order-ideal in $P_{k,Z}(m,n)$, the M\"obius
function on $P_{k,Z'}(m,n)$ is the restriction of the M\"obius
function on $P_{k,Z}(m,n)$ to $P_{k,Z'}(m,n)$. Now the statement
follows from Lemma \ref{L:inducedposet}.
\end{pf}

We now prove the formula of Conjecture \ref{C:mainmoeb} under the
condition that $Z$ contains a $k-1$-assignment.

\begin{Thm}\label{T:k-1}
If $Z$ contains a $k-1$-assignment, then
\begin{equation} \label{eq:mainthm} 
F_{k,Z}(m,n)=\sum_{\al\in B(m,n)}{\frac{g(\al)}{(m-i_\al)(n-j_\al)}}.
\end{equation}
\end{Thm}

\begin{pf}
If $Z$ contains a $k$-assignment, the statement is true. We will
therefore assume that $Z$ does not contain a $k$-assignment. We
will use induction on $k$, and for a specific value of $k$, we
will also use induction on the number of elements \emph{not}
belonging to $Z$. Hence we will assume that the statement has been
proved for $F_{k,Z'}$, for every superset $Z'$ of $Z$.

Notice that if there is a row $r$ which belongs to every
$k-1$-covering of $Z$, then by Theorem \ref{T:deleterow}, $$
F_{k,Z}(m,n)=F_{k-1,Z\setminus r}(m-1,n).$$ By induction on $k$,
this is $$ \sum_{\al\in P_{k-1,Z\setminus
r}(m-1,n)\setminus\hat{0}
}{\frac{g(\al)}{(m-1-i_\al)(n-j_\al)}}.$$ But $P_{k-1,Z\setminus
r}(m-1,n)$ is isomorphic to $P_{k,Z}(m,n)$ via the isomorphism
$\al\mapsto \al\cup \{r\}$, and $i_{\al\cup \{r\}}=i_\al+1$. Hence
$F_{k-1,Z\setminus r}(m-1,n)$ equals \eqref{eq:mainthm}.

We can assume that the $k-1$-assignment in $Z$ is contained in the
first $k-1$ rows and the first $k-1$ columns. Actually, we can
assume that all elements of $Z$ are contained in this $k-1$ by
$k-1$ square, since if there is an element of $Z$ in another
column, say, then the row containing this element must be used in
every $k-1$-covering of $Z$, and we can use induction on $k$.

We now subtract the minimum of the elements in columns
$k,\dots,n$, and sum over the $m$ possibilities for the row $r$
containing the new zero, using induction:

$$
F_{k,Z}(m,n)=\frac{1}{m(n-k+1)}+\frac{1}{m}\sum_{r}{\sum_{\begin{Sb}
\al\in P_{k,Z}(m,n)\setminus\hat{0}\\ r\in \al \end{Sb}
}{\frac{g(\al)}{(m-i_\al)(n-j_\al)}}}.$$

Now we change the order of summation, and use the fact that the
number of rows $r$ for which $r\in \al$ is $i_\al$.

\begin{multline}
F_{k,Z}(m,n)=\frac{1}{m(n-k+1)}+\frac{1}{m}\sum_{\al\in
P_{k,Z}(m,n)\setminus\hat{0}}{\frac{i_\al g(\al)}{(m-i_\al)(n-j_\al)}}=\\
=\frac{1}{m(n-k+1)}+\sum_{\al\in
P_{k,Z}(m,n)\setminus\hat{0}}{\frac{g(\al)}{(m-i_\al)(n-j_\al)}}-\sum_{\al\in
P_{k,Z}(m,n)\setminus\hat{0}}{\frac{g(\al)}{m(n-j_\al)}}.
\end{multline}

We have to show that the first term and the last sum cancel, that
is, that $$ \sum_{\al\in
P_{k,Z}(m,n)\setminus\hat{0}}{\frac{g(\al)}{(n-j_\al)}}=\frac{1}{n-k+1}.$$

When $\al$ is the set of the first $k-1$ columns, we get a term
which equals the right hand side. Hence we are going to show that
everything else cancels. It clearly suffices to show that if $J$
is a proper subset of the set of the first $k-1$ columns, then

$$\sum_{\begin{Sb} \al\in P_{k,Z}(m,n)\\ column-set(\al)=J
\end{Sb} }{g(\al)}=0.$$

If $J$ is nonempty, then the induced sub-poset of $P_{k,Z}(m,n)$
consisting of those elements that contain $J$ is isomorphic to
$P_{k-\left|J\right|, Z\setminus J}(m,n-\left|J\right|)$ via the
isomorphism $\al\mapsto \al\setminus J$. Hence the claim follows by
induction on $k$.

If on the other hand $J=\O$, we have

\begin{multline}
\sum_{\begin{Sb} \al\in P_{k,Z}(m,n)\\
column-set(\al)=\O\end{Sb}}{g(\al)} =\sum_{\al\in
P_{k,Z}(m,n)\setminus\hat{0}}{g(\al)}-\sum_{\begin{Sb} \al\in
P_{k,Z}(m,n)\setminus\hat{0}\\ column-set(\al)\neq\O
\end{Sb}}{g(\al)}.
\end{multline}

By what we already know, the second sum equals 1, since everything
cancels except the term where $\al$ is the set of the first $k-1$
columns. If we let $\al_0$ be the intersection of all coverings of
$Z$ with $k-1$ rows and columns, then $$\sum_{\al\in
P_{k,Z}(m,n)\setminus\hat{0}}{g(\al)}=\sum_{\al\leq
\al_0}{g(\al)}=1,$$ by the definition of $g$. Hence
$$\sum_{\begin{Sb} \al\in P_{k,Z}(m,n)\\ column-set(\al)=\O
\end{Sb} }{g(\al)}=1-1=0,$$ which proves the theorem.
\end{pf}

As a special case Theorem \ref{T:k-1} enables us to compute the
expected cost for completing a zero cost $k-1$-assignment, thus
answering a question posed in \cite{CS}

\begin{Cor}\label{C:diagonal}
Let $Z_k=\{(1,1),(2,2),\dots,(k-1,k-1)\}$. Then
\[
F_{k,Z_k}(m,n)=\sum_{i+j<k}{\binom{k-1}{i,j,k-1-i-j}
\frac{(-1)^{k-1-i-j}}{(m-i)(n-j)}}.
\]
In particular,
\[
F_{k,Z_k}(k,n)=\frac{1}{k}\sum_{j=0}^{k-1}{\frac{1}{n-j}},
\]
and
\[F_{k,Z_k}(k,k)=\frac{1}{k}\left(1+1/2+\dots+1/k\right)\sim\frac{\log
k}{k}.\]
\end{Cor}

\begin{pf}
The $k-1$-coverings of $Z_k$ are all obtained by choosing, for
each zero, whether to cover it with a row or with a column. It is
easily seen that all subsets of $k-1$-coverings (that is, all
partial $k-1$-coverings) can be obtained as intersections of
$k-1$-coverings. The elements of the poset $P_{k,Z_k}$ are
therefore obtained by choosing, for every zero, whether to cover
it with a row, with a column, or not to cover it. For any partial
$k-1$-coverings $\hat 0 \neq \al\le \be$, the interval
$[\al,\be]\subset P_{k,Z_k}$ is isomorphic to the Boolean lattice
on $|\al|-|\be|$ elements. We get
\begin{multline}
\mu(\hat 0,\hat 1)=-\sum_{\al > \hat 0}\mu(\al,\hat 1) =-\sum_{\al
> \hat 0}(-1)^{|\al|}=\\ =-\sum_{l=0}^{k-1} (-1)^l
\binom{k-1}{l}2^l=-(1-2)^{k-1}=(-1)^k.
\end{multline}

Any interval $[\hat 0, \al]\subset P_{k,Z_k}$, $\al\neq \hat 0$ is
isomorphic to $P_{l,Z_l}$, where $l=k-|\al|$. Hence $\mu(\hat
0,\al)=(-1)^{k-|\al|}$, for all $\al \in P_{k,Z_k}.$ The number of
partial $k-1$-coverings with $i$ rows and $j$ columns is given by
the trinomial coefficient $\binom{k-1}{i,j,k-1-i-j}$. This proves
the first statement.

For the second statement, observe that we can rewrite $F_{k,Z_k}(m,n)$
by
\begin{multline}
\sum_{j=0}^{k-1}\sum_{i=0}^{k-1-j}{\binom{k-1}{i,j,k-1-i-j}
\frac{(-1)^{k-1-i-j}}{(m-i)(n-j)}}=\\
\sum_{j=0}^{k-1}\binom{k-1}{j}\frac{1}{n-j}\sum_{i=0}^{k-1-j}
{\binom{k-1-j}{i}\frac{(-1)^{k-1-i-j}}{m-i}}=\\
[\text{by Lemma \ref{L:sum} }]
=\sum_{j=0}^{k-1}\frac{\binom{k-1}{j}}
{\binom{m-1}{k-1-j}}\cdot\frac{1}{m(n-j)}.
\end{multline}

In particular, when $k=m$, we have
\[
F_{k,Z_k}(k,n)=\frac{1}{k}\sum_{j=0}^{k-1}{\frac{1}{n-j}}.
\]
\end{pf}

\section{Some consequences of the conjecture}\label{S:consequences}
A consequence of the Main Conjecture is the following.

\begin{Thm} \label{T:ettor}
Let $k\le m,n$ be positive integers and $Z$ the set of zero
entries. If Conjecture \ref{C:main} is true, then there exist
integers $b_{i,j}$ (independent of $m$ and $n$) such that

\begin{equation} 
\fkz=\sum_{i+j<k}\frac{b_{i,j}}{(m-i)(n-j)},
\end{equation}

Moreover, if all the zeros in $Z$ are in an $m' \times n'$ sub-array $A'$,
then $b_{i,j}=1$ if $i\ge m'$ or $j\ge n'$.
\end{Thm}

\begin{pf}
Assume Conjecture \ref{C:main}.
Rewrite $\fkz$ as $$\sum_{x+y<k}{\frac{1-q_{x,y}}{(m-x)(n-y)}},$$
where $q_{x,y}$ is the probability that a random set with $x$ rows
and $y$ columns (chosen uniformly) is a bad set. Now refine the
counting of bad sets by letting $s$ and $t$ be the number of rows
resp. columns that intersect $A'=[1,m']\times [1,n']$. Then
\begin{equation}\label{eq:qxy}
\frac{q_{x,y}}{(m-x)(n-y)}=\frac{\sum_{s=0}^x \sum_{t=0}^y
{\binom{m-m'}{x-s}\binom{n-n'}{y-t}d_{s,t}^{k-1-y-x}(A')}}
{\binom{m}{x}(m-x)\binom{n}{y}(n-y)},
\end{equation}
where $d_{s,t}^r(A')$ is the number of sets of $s$
rows and $t$ columns that are not partial $s+t+r$-coverings.

We use partial fractions with respect to $m$ to get
\begin{equation}\label{eq:parfrac}
\frac{\binom{m-m'}{x-s}}{\binom{m}{x}(m-x)}=
\sum_{i=0}^{\min\{x,m'-1\}}
\frac{(-1)^{s-i}\binom{x}{i}\binom{m'-i+x-s-1}{m'-i-1}}{m-i}.
\end{equation}

Inserting \eqref{eq:parfrac} and the
corresponding identity for
$\frac{\binom{n-n'}{y-t}}{\binom{n}{y}(n-y)}$
into \eqref{eq:qxy} we get that the contribution of
$\frac{q_{x,y}}{(m-x)(n-y)}$ to $b_{i,j}$ is the integer
\[
(-1)^{i+j}
\tbinom{x}{i}\tbinom{y}{j} \sum_{s=0}^{x}
\sum_{t=0}^{y}
(-1)^{s+t}\tbinom{m'-i+x-s-1}{m'-i-1}\tbinom{n'-j+y-t-1}{n'-j-1}
d_{s,t}^{k-1-y-x}(A'),
\]
if $i\le \min \{x,m'-1\}$ and $j\le \min\{y,n'-1\}$. Otherwise
the contribution is zero and the theorem follows.
\end{pf}

Theorem \ref{T:ettor} enables us to write the value of $F_{k,Z}$ as a
triangle of $b_{i,j}$'s which has been very convenient and is the notation
used in the Appendix. For example the resulting expected value in
Example \ref{E:poset} can be written as

\[
\begin{picture}(924,924)(1789,-1273)
\thinlines
\put(1801,-661){\framebox(900,300){}}
\put(1801,-1261){\framebox(300,900){}}
\put(1801,-961){\framebox(600,600){}}
\put(2551,-586){\makebox(0,0)[lb]{\smash{\SetFigFont{6}{7.2}
{\rmdefault}{\mddefault}{\updefault}1}}}
\put(1951,-1186){\makebox(0,0)[lb]{\smash{\SetFigFont{6}{7.2}
{\rmdefault}{\mddefault}{\updefault}1}}}
\put(1951,-886){\makebox(0,0)[lb]{\smash{\SetFigFont{6}{7.2}
{\rmdefault}{\mddefault}{\updefault}1}}}
\put(1846,-586){\makebox(0,0)[lb]{\smash{\SetFigFont{6}{7.2}
{\rmdefault}{\mddefault}{\updefault}-2}}}
\put(2236,-586){\makebox(0,0)[lb]{\smash{\SetFigFont{6}{7.2}
{\rmdefault}{\mddefault}{\updefault}1}}}
\put(2221,-886){\makebox(0,0)[lb]{\smash{\SetFigFont{6}{7.2}
{\rmdefault}{\mddefault}{\updefault}1}}}
\end{picture}.
\]
\medskip

The equation in the preceding proof can be used to give exact
formulas for $b_{i,j}$. Let $r=x+y$ and sum over all $x,y$, then
\begin{multline}\label{eq:bij}
b_{i,j}=
1-(-1)^{i+j}\sum_{r=i+j}^{k-1} \sum_{s=0}^{m'-k+r}
\sum_{t=0}^{n'-k+r} \bigg( (-1)^{s+t} d_{s,t}^{k-1-r}(A')\cdot\\
\cdot\sum_{x=\max\{i,s\}}^{r-\max\{j,t\}} \tbinom{x}{i}\tbinom{r-x}{j}
\tbinom{m'-i+x-s-1}{m'-i-1}\tbinom{n'-j+r-x-t-1}{n'-j-1}\bigg ).
\end{multline}
Here we have used the fact that $d_{s,t}^{k-1-r}(A')=0$ if
$s+k-1-r\ge m'$ or $t+k-1-r\ge n'$.

\subsection{The probability that the smallest element in a row is used}

In \cite{Olin}, B.~Olin conjectured, on the basis of experimental
data, that in the case $k=m=n$, the probability that the smallest
element in a row is used in the optimal $n$-assignment tends to
1/2, as $n\to \infty$ (Olin was considering uniform distribution
in $[0,1]$ instead of exponential distribution, but this does not
seem to matter in the limit). We now show that Conjecture
\ref{C:main} implies the conjecture of Olin.

\medskip

Note that the probability that the smallest element in a row is
used is equal to the probability that the zero is used in a matrix
with just one zero. This in turn is equal to the probability that
the smallest element in a matrix without zeros is used. Assuming
Conjecture \ref{C:main}, we can compute this probability.

\begin{Thm}
Let $A$ be an $m\times n$-matrix with one zero. Let $k$ be a
positive integer. Assuming Conjecture \ref{C:main}, the
probability that the zero is used in the optimal $k$-assignment is
$$ 1-\frac{\binom{k}{2}}{mn}.$$
\end{Thm}

\begin{pf}
We can assume that the zero is $a_{1,1}$, so $Z=\{(1,1)\}$.
Note that one element
from the first row is always used in the optimal $k$-assignment.
We therefore compute the probability that another element in the
first row, say $a_{1,2}$, is used. Let $Z'=\{(1,1),(1,2)\}$. Let
$p_{i,j}$ and $p'_{i,j}$ be the corresponding probabilities for
partial $k-1$-coverings with $i$ rows and $j$ columns. We are interested in
the difference $p_{i,j}-p'_{i,j}$. If $i+j<k-1$, then
$p_{i,j}=p'_{i,j}=1$. If $i+j=k-1$, then $p_{i,j}-p'_{i,j}$ is the
probability that, if $i$ rows and $j$ columns are chosen randomly,
position $(1,1)$ is covered, but position $(1,2)$ is not. This is the
case if and only if column 1 is used, but not row 1 or column 2.
The probability for this is
$$\frac{m-i}{m}\cdot\frac{n-j}{n}\cdot\frac{j}{n-1}=\frac{j(m-i)(n-j)}{
mn(n-1)}.$$

We have
\[
Pr[\text{$a_{1,2}$ is
used}]=\sum_{i+j<k}{\frac{p_{i,j}-p'_{i,j}}{(m-i)(n-j)}}=\sum_{i+j=
k-1}{\frac{j}{mn(n-1)}}=
\frac{\binom{k}{2}}{mn(n-1)}.
\]

Hence the probability that one of the elements $a_{1,2},
a_{1,3},\dots, a_{1,n}$ is used, is $$\frac{\binom{k}{2}}{mn}.$$

\end{pf}

\section{Some more evidence for the Main Conjecture}\label{s:evidence}
If Conjecture \ref{C:main} is true then it should of course be
consistent with the recursion from Theorem \ref{T:recursion} in the 
case there is no doubly covered entry.
In this
section we will prove consistency of $b_{i,j}$ for the special case
when $i+j=k-1$. In this case the formula above simplifies to

\[
b_{i,j}=1-(-1)^{k-1} \sum_{s=0}^{m'-1} \sum_{t=0}^{n'-1}
(-1)^{s+t} d_{s,t}^{0}(A')
\binom{m'-s-1}{m'-i-1}\binom{n'-t-1}{n'-j-1}.
\]

One can of course do the same calculations starting with $p_{i,j}$
instead of $q_{i,j}$ and get that $b_{i,j}$ is

\begin{equation}\label{eq:pij}
b_{i,j}=(-1)^{k-1} \sum_{s=0}^{m'-1} \sum_{t=0}^{n'-1} (-1)^{s+t}
g_{s,t}^{0}(A') \binom{m'-s-1}{m'-i-1}\binom{n'-t-1}{n'-j-1},
\end{equation}
where $g_{s,t}^{r}(A')$ is the number of partial $s+r+t$-coverings of
$A'$, with $s$ rows and $t$ columns. This is the formula we will use.

We do the case when $Z\subset [1,m']\times [1,n']$ is optimally
covered by $n'$
columns. Let $B_{x}$, $x=1\ldots m'$ be the matrix with zeros in
$Z\cup {(x,n'+1)}$. Let also $C$ be the
matrix with zeros in $Z\cup {(m'+1,n'+1)}$.
Then the recursion step corresponding to the optimal covering of
columns gives

\begin{equation}\label{eq:rec}
b_{i,j}(A)=\frac{\sum_{x=1}^{m'} b_{i,j}(B_{x})}{i}+\frac{(i-m')
b_{i,j}(C)}{i}= \frac{\sum_{x=1}^{m'}
b_{i,j}(B_{x})-m'b_{i,j}(C)}{i}+b_{i,j}(C).
\end{equation}

Note that we can include an extra row and column in $A'$ if we just
replace $m'$ by $m'+1$ and $n'$ by $n'+1$ in equation \eqref{eq:pij}.
To simplify the
calculations we will do this and thereby assume that $A', B_x'$ and
$C'$ all are $m'+1\times n'+1$ submatrices of $A,B_x$ and $C$
respectively.

We will now prove that equation \eqref{eq:pij} remains valid
during a step of the recursion. Partition
$g_{s,t}^{0}(X)=g_{s,t}^{0,c}(X)+g_{s,t}^{0,r}(X)+g_{s,t}^{0,0}(X)$,
by

$g_{s,t}^{0,c}(X)=$ the number of partial $s+t$-coverings of $X$ that use
the last column.

\smallskip
$g_{s,t}^{0,r}(X)=$ \parbox{4in}{the number of partial
$s+t$-coverings of $X$ that do
not include the last column but the last row.}

\smallskip
$g_{s,t}^{0,0}(X)=$ \parbox{4in}{the number of partial
$s+t$-coverings of $X$ that include neither the last row nor the last
column.}

\smallskip
We can now establish a number of equalities from which the result
will follow. First note $g_{s,t}^{0,0}(C')=0$. Second
$g_{s,t}^{0,c}(A')+g_{s,t}^{0,r}(A')=
g_{s,t}^{0,c}(C')+g_{s,t}^{0,r}(C')$,
this is used to match away the last term in \eqref{eq:rec} and the
``non-wanted'' coverings of the enlarged $A'$.

Third we have
$\sum_{x=1}^{m'}g_{s,t}^{0,c}(B_{x}')=m'g_{s,t}^{0,c}(C'))$.

\smallskip
We now put together

$\sum_{x=1}^{m'}g_{s+1,t}^{0,r}(B_{x}')=\sum_{x=1}^{m'}g_{s,t}^{0,0}(B_{x}')
=sg_{s,t}^{0,0}(A'))$ and

$g_{s+1,t}^{0,r}(C')=g_{s,t}^{0,0}(A'))$
to get the last identity needed

\[
\frac{\sum_{x=1}^{m'} \binom{m'-s}{m'-i} g_{s,t}^{0,0}(B_{x}')-
\binom{m'-(s+1)}{m'-i} (\sum_{x=1}^{m'}
g_{s+1,t}^{0,r}(B_{x})-m'g_{s+1,t}^{0,r}(C))}{i}=
\]
\[
\binom{m'-s}{m'-i}g_{s,t}^{0,0}(A').
\]

Plugging these identities and \eqref{eq:pij} into
recursion \eqref{eq:rec} everything will match nicely.

\section{Asymptotic consequences}\label{S:asymptotics}

In this section, we show that certain limits can be computed from
Conjecture \ref{C:main}. We assume throughout this section that
Conjecture \ref{C:main} holds.

Let $K$ be a subset of the unit square $[0,1]\times [0,1]$. Let
$n$ be a positive integer. Later, we will let $n$ tend to
infinity.

We divide the unit square into squares $s_{i,j}$ of side $1/n$,
such that
\[
s_{i,j}=\left[\frac{i-1}{n},\frac{i}{n}\right]\times
\left[\frac{j-1}{n},\frac{j}{n}\right].
\]

Let $Z_n$ be the set of all $(i,j)$ such that $s_{i,j}$ intersects
$K$.

We shall show that under certain circumstances, we can compute the
limit of $F_n=F_{n,Z_n}(n,n)$, as $n\to \infty$.

For a fixed $n$ we have
\[
F_n=\sum_{i,j}{\frac{p_{i,j}}{(n-i)(n-j)}},
\]
where $p_{i,j}$ is the probability that a set of $i$ rows and $j$
columns is a partial $k-1$-covering.

By the substitutions $x=i/n$, $y=j/n$, we can rewrite this as
\[
\int_{0}^{1}{\int_{0}^{1}{\frac{p_n(x,y)}{(1-x)(1-y)}\, dx\,dy}},
\]
where $p_n(x,y)=p_{\lfloor xn\rfloor,\lfloor yn\rfloor}$.

If $S$ is a subset of ${\mathbb{R}}^2$, we let $M_{a,b}(S)$ be the
infimum of $a\lambda(S_1)+b\lambda(S_2)$ (Lebesgue-measure), taken
over all sets $S_1, S_2$ of real numbers such that $S\subseteq
\left(S_1\times {\mathbb{R}}\right) \cup \left({\mathbb{R}}\times
S_2\right)$.

\begin{Lemma}
Suppose that $K$ is compact. Let $x$ and $y$ be such that
$M_{1-x,1-y}(K)<1-x-y$. Then $p_n(x,y)\to 1$, as $n\to\infty$.
\end{Lemma}

\begin{pf}
We can find $S_1$, $S_2$ and $\delta$ such that $K\subseteq
\left(S_1\times {\mathbb{R}}\right) \cup \left({\mathbb{R}}\times
S_2\right)$, and
$(1-x)\lambda(S_1)+(1-y)\lambda(S_2)<1-x-y-\delta$. By the
compactness of $K$, we can assume that $S_1$ and $S_2$ are finite
unions of intervals.

Let $n$ be a (large) integer. Then the number of intervals of the
form $\left[\frac{i-1}{n}, \frac{i}{n}\right]$ which intersect
$S_1$ is at most $n\cdot\lambda(S_1)+constant$, and similarly for
the number of such intervals which intersect $S_2$. Hence we can
find a covering $\al$ of $Z_n$ using at most
$n\cdot(\lambda(S_1)+\lambda(S_2))+constant$ rows and columns.

Now choose a set $\al'$ consisting of $\lfloor nx \rfloor$ rows and
$\lfloor ny \rfloor$ columns at random. The average number of rows
and columns in $\al$ which are not in $\al'$ is at most
$$n\left((1-x)\lambda(S_1)+(1-y)\lambda(S_2)\right)+constant <
n(1-x-y-\delta) + constant.$$ If $n$ is large, then by the weak
law of large numbers, the probability that $\left|\al-\al'\right|$
deviates from its mean value by more than $n\cdot\delta$ tends to
zero. With high probability, $\left|\al \bigcup \al'\right| < n$.
Hence with high probability, $\al'$ is a partial $n-1$-covering.
\end{pf}

We say that a covering of $Z_n$ is \emph{minimal}, if no proper
subset is a covering of $Z_n$.

\begin{Lemma}
Suppose that the number of minimal coverings of $Z_n$ grows slower
than $(1+\epsilon)^n$ for every $\epsilon$, as $n\to \infty$. Let
$x$ and $y$ be such that $M_{1-x,1-y}(K)>1-x-y$. Then $p_n(x,y)\to
0$, as $k\to\infty$.
\end{Lemma}

\begin{pf}
Let $T(n)$ be the number of minimal coverings of $Z_n$. Let $\al$ be
a minimal covering. Since $M_{1-x,1-y}(K)>1-x-y+\delta$, for a
suitably chosen $\delta>0$, we have
$$(1-x)i_\al+(1-y)j_\al>n\cdot(1-x-y+\delta).$$ If we let $i=\lfloor
n\cdot x\rfloor$ and $j=\lfloor n\cdot y\rfloor$, we certainly
have $$(1-x)\cdot i_\al+(1-y)\cdot j_\al>n-i-j+2\delta_1\cdot n,$$ for
some $\delta_1$ slightly smaller than $\delta/2$.

Let $\al'$ be a random set constructed by letting each row belong to
$\al'$ with probability $x$, and each column with probability $y$.
Since with probability at least 1/4, $\al'$ will contain at least
$i$ rows and $j$ columns, in order to prove that $p_n(x,y)$ tends
to 0, it will suffice to show that the probability that one can
cover $Z_n$ by adding $n-i-j-1$ rows and columns to $\al'$ tends to
0.

Let $B(p)$ denote the random variable that is $1$ with probability $p$ and
$0$ with probability $1-p$.
The number of rows in $\al\setminus \al'$ is a sum of $i_\al$
independent $B(1-x)$-distributed random variables, and similarly,
the number of columns in $\al\setminus \al'$ is a sum of $j_\al$
independent $B(1-y)$-variables. In order that the size of
$\al\setminus \al'$ be at most $n-i-j-1$, one of these sums must
deviate by at least $n\cdot \delta_1$ from its mean value. By
standard estimates in probability theory, this is at most
$e^{-c\cdot n}$, for some constant $c$ depending only on
$\delta_1$.

Hence the probability that one can obtain a covering of $Z_n$ by
adding at most $n-i-j-1$ rows and columns to $\al'$ is at most
$T(n)\cdot e^{-c\cdot n}$, which, by assumption, tends to 0 as
$n\to \infty$.
\end{pf}

\begin{Cor} \label{Cor:lim}
If $K$ is compact, and the number of minimal coverings of $Z_n$
grows slower than $(1+\epsilon)^n$ for every $\epsilon$, as $n\to
\infty$, then
\[
\lim_{n\to\infty}{F_n}=\int_D{\frac{dx\,dy}{(1-x)(1-y)}},
\]
where $D$ is the region in which $M_{1-x,1-y}(K)<1-x-y$.
\end{Cor}

\subsection{An example: Zeros outside a quarter of a circle}
Let $K$ be the region $x^2+y^2\geq 1$,
$0\leq x,y\leq 1$. To find
$M_{a,b}(K)$, for given $a$ and $b$, we have to find the point
where the slope of the curve $x^2+y^2=1$ is $-a/b$. The part of
the curve which is above this point should be covered with a
horizontal strip, and the part which is to the right of this point
should be covered with a vertical strip. This point is
\[
\left(\frac{a}{\sqrt{a^2+b^2}},\frac{b}{\sqrt{a^2+b^2}}\right).
\]
Hence
\[
M_{a,b}(K) =
a\left(1-\frac{a}{\sqrt{a^2+b^2}}\right)+b\left(1-\frac{b}{\sqrt{a^2+b^2}
}\right
),
\]
which simplifies to
\[
a+b-\sqrt{a^2+b^2}.
\]

The number $T(n)$ of minimal coverings of $Z_n$ is at most $n+1$,
since any minimal covering consists of, for some $i$, $0\leq i\leq
n$, columns $i,\dots, n$, and all rows that contain some element
of $Z_n$ not covered by the columns.

The region $D$ is given by the inequality $M(K_{1-x,1-y})<1-x-y$,
which becomes
\[
(1-x)^2+(1-y)^2>1.
\]
In other words, $D$ is the part of the unit square which is
outside the circle of radius 1 centered in the point $(1,1)$. By
Corollary \ref{Cor:lim}, we have
\[
\lim_{n\to\infty}{F_n}=\int_{D}{\frac{dx\,dy}{(1-x)(1-y)}}=\frac{\pi^
2}{24}
.
\]

\subsection{Generalizing the example} Fix $p>1$. Let $K$ be the region
given by $x^p+y^p\geq 1$, $0\leq x,y\leq 1$. We wish to compute
$M_{a,b}(K)$. To find the point on the curve $x^p+y^p=1$ where the
slope is $-a/b$, we compute:
\[
\frac{dy}{dx}=-x^{p-1}(1-x^p)^{1/p-1}=-(x/y)^{p-1}.
\]
It is not difficult to see that the point is
\[
\left (\frac{a^{\frac{1}{p-1}}}{\left(a^{\frac{p}{p-1}}+
b^{\frac{p}{p-1}}\right)^{1/p}}, \frac{b^{\frac{1}{p-1}}}
{\left(a^{\frac{p}{p-1}}+b^{\frac{p}{p-1}}\right)^{1/p}}
\right ).
\]
We have
\begin{multline}
M_{a,b}(K)=a\left(1-\frac{a^{\frac{1}{p-1}}}{\left(a^{\frac{p}{p-1}}+
b^{\frac{p}{p-1}}\right)^{1/p}}\right)+
b\left(1-\frac{b^{\frac{1}{p-1}}}{\left(a^{\frac{p}{p-1}}+
b^{\frac{p}{p-1}}\right)^{1/p}}\right)=\\
=a-\frac{a^{\frac{p}{p-1}}}{\left(a^{\frac{p}{p-1}}+
b^{\frac{p}{p-1}}\right)^{1/p}}+b-\frac{b^{\frac{p}{p-1}}}{\left(a^{\frac
{p}{p-1
}}+ b^{\frac{p}{p-1}}\right)^{1/p}}=\\
=a+b-\left(a^{\frac{p}{p-1}}+
b^{\frac{p}{p-1}}\right)^{\frac{p-1}{p}}.
\end{multline}

Hence
the inequality $M_{1-x,1-y}(K)< 1-x-y$ becomes
\[
\left((1-x)^{\frac{p}{p-1}}+
(1-y)^{\frac{p}{p-1}}\right)^{\frac{p-1}{p}}> 1.
\]
Again the number of minimal coverings of $Z_n$ is at most $n+1$.
To find the limit of $F_n$, we now by Corollary \ref{Cor:lim}
wish to compute the integral
\[
\int_{x=0}^{1}{\int_{y=(1-x^{u})^{1/u}}^{1}{\frac{dx\,dy}{xy}}},
\]
where $u=p/(p-1)$. We eliminate the inner integral, and obtain
\[
\int_0^1{\frac{-\log\left((1-x^u)^{1/u}\right)}{x}\,dx}=
-\frac{1}{u}\int_0^1{\frac{\log(1-x^u)}{x}\,dx}.
\]

We make the substitution $t=x^u$, which
gives $dx = \frac{dt}{ux^{u-1}}$. We get
\begin{multline}
-\frac{1}{u^2}\int_0^1{\frac{\log(1-t)}{t}\,dt}
=-\frac{1}{u^2}\int_0^1{\left(\frac{-t-t^2/2-t^3/3-\dots}{t}\right)\,dt}
=\\
=\frac{1}{u^2}\left(1+\frac{1}{4}+\frac{1}{9}+\dots\right)
=\frac{1}{u^2}\cdot\frac{\pi^2}{6}.
\end{multline}
Substituting back $p/(p-1)$ for $u$,
\[
\lim_{n\to\infty}{F_n}=\left(1-\frac{1}{p}\right)^2\cdot\frac{\pi^2}{6}
.
\]

\section{Remarks}

Looking at the triangles of $b_{i,j}$'s one discover certain patterns.
Define the property {\it
acyclic} recursively by first letting the empty set be acyclic.
Second, let $Z$ be an acyclic set
that has an optimal covering $\al$ with only rows or only columns.
Let $(i,j)$ be any position not covered by $\al$. Then we define
$Z\cup (i,j)$ to also be acyclic.

\begin{Cor}\label{C:b00}
Assume Conjecture \ref{C:main}. If $Z$ is acyclic,
then $b_{0,0}=(-1)^{|Z|}\binom{k-1}{|Z|}$. \qed
\end{Cor}

This corollary can be proved by first specializing \eqref {eq:bij} to
$b_{0,0}$ and then study the set of bad sets and find the proper
bijections which proves it by induction on $|Z|$.
The proof we have
found is much in the spirit of the proof in Section \ref{s:evidence}
of consistency for $b_{i,j}$, $i+j=k-1$.
There are some highly non-trivial details and we omit the proof.

\medskip
Given a set of zeros $Z$, let $Z_{row}(i)$, be the number of zeros in row
$i$ and $Z_{col}(j)$ the number of zeros in column $j$.
Ordering the nonzero $Z_{row}(i)$, $1\le i\le m$ by size, we get an
integer partition of $|Z|$, call it $\lambda_{row}(Z)$. Similarly we
define $\lambda_{col}(Z)$.

\begin{Conjecture}
Given $k$ and $Z$, write $\fkz=\sum_{i+j<k}\frac{b_{i,j}}{(m-i)(n-j)}$
for some (uniquely defined) integers $b_{i,j}$
as in Theorem \ref{T:ettor}.
Then the value of $b_{0,0},b_{1,0},\dots,b_{k-1,0}$ only depend on
$\lambda_{row}(Z)$ and $b_{0,0},b_{0,1},\dots,b_{0,k-1}$ only depend on
$\lambda_{col}(Z)$.
\end{Conjecture}
It is not difficult to see that this conjecture
implies Corollary \ref{C:b00}.

\medskip
We mention that Theorem \ref{T:e-d} can be generalized in the
following way:

\begin{Thm}
Let $a_i$ and $b_i$, $1\leq i\leq N$, be real numbers. Let $X$ be
an exponentially distributed random variable with mean 1. Define
the random variable $I\in \{1,\dots,N\}$ by $a_I+b_IX=\min_i
{a_i+b_iX}$. Then $\text{E}(b_I)=\text{E}(a_I+b_IX)-\min(a_i)$.
\end{Thm}

\begin{pf}
Let $f(x)=\min_{i}(a_i+b_ix)$. Then
\begin{multline}
-\min_i(a_i)=-f(0)=\int_0^\infty{\frac{d}{dx}\left(e^{-x}f(x)\right)
\,dx}=\\
=\int_0^\infty{e^{-x}f'(x)\,dx}-\int_0^\infty{e^{-x}f(x)\,dx}=\\
=\text{E}(f'(X))-\text{E}(f(X))=\text{E}(b_I)-\text{E}(a_I+b_IX).
\end{multline}
\end{pf}

As a last remark we want to point at the recently published preprint 
\cite{BCR00}, where they generalize the random assignment problem in a 
different direction.  They assume that the matrix has independent 
exponential random variables with different intensities.  In the case 
when the intensities form a rank 1 matrix they give a conjecture of 
the expected value of the optimal $k$-assignment.  Observing that the 
binomial coefficient in their conjecture is the M\"obius function of a 
truncated Boolean algebra, it is not difficult to join our conjecture 
with theirs to the case with a rank 1 matrix of intensities, but in 
some entries the random variables are replaced by zeros.  We have done 
no computations to see if this join of conjectures is true.

\newpage

\section{Appendix: Computer generated verification of Conjecture
\ref{C:CS} for $k\leq 5$}

Our algorithm for calculating $F_{k,Z}(m,n)$ is essentially a
formalization of the method outlined in Section
\ref{S:computation}. It was implemented in Maple to take advantage
of Maple's built in functions for simplifying rational
expressions.

Let us immediately state that we do \emph{not} have a proof that
this algorithm always returns a result. However, if a result is
returned, we know that this result is correct. The difficulty lies
in that this algorithm does not use fixed values for $m$ and $n$,
but tries to compute $F_k(m,n)$ as a rational function of $m$ and
$n$.  For fixed numerical values of $m$, $n$ and $k$, it is in
principle always possible to calculate $F_k(m,n)$.

The algorithm uses a data-structure which we will here call an
\emph{array}. An array is a matrix whose elements are finite sums
of the form $$\sum_i{r_iX_i},$$ where $X_1,X_2,X_3,\dots$ are
independent exp(1)-distributed random variables, and $r_i$ is a
rational function in the variables $m$ and $n$. An array $M$ is
thought of as representing a window of an $m$ by $n$ matrix $A$,
all of whose entries outside the region visible in $M$ are
exp(1)-variables, independent among themselves, and independent of
all entries in $M$. In this way we can for example represent a
pattern $Z$ of zeros, independently of $m$ and $n$, and compute
$F_{k,Z}(m,n)$ as a rational function in $m$ and $n$.

Below are the cases computed to calculate $F_5(m,n)$ using Theorem
\ref{T:recursion}. The same computer program has been used to
verify Conjecture \ref{C:CS} also for $k=6$, but this computation
needed about 10,000 cases, and is not written out here. In each
case, $M$ is an array, and an optimal covering of the zeros of $M$
is computed. The other cases that are used for a particular case
are listed. The expected value of the optimal $k$-assignment is
denoted $F$. The value of $F$ is given as a matrix. This matrix
represents the coefficients $b_{i,j}$. Hence the entry in position
$(i+1,j+1)$ of the matrix is the coefficient of $1/((m-i)(n-j))$
in the expression for $F$. Sometimes the value of $F$ cannot be
written as a linear combination of such terms. Then the extra
terms are written out explicitly.

\subsection{The case $k=m=n=7$}

In case $k=m=n$, we can take advantage of the fact that every
$k$-assignment contains one element from every row and every
column. This means that we can start by subtracting the minimal
element from every row, and then continue by subtracting the
minimal element from every column that does not already contain a
zero. Hence we need only take into account the cases where $Z$
contains at least one element in every row and every column. It
turns out that in the calculation of $F_7(7,7)$, all the cases
that occur can easily be reduced to cases that have already
occurred in the calculation of $F_k(m,n)$ for $k\le 5$. In this way
we have verified that $F_7(7,7)=1+1/4+1/9+\dots+1/49$, as
conjectured.

\newpage

\tiny
\input {Out2.tex}

\end{document}

%% file: Out2.tex
\begin{maplettyout}
\begin{center}{\bf Case   1}\end{center}
\end{maplettyout}

\begin{maplelatex}
\[
k=0, \,M=[], \,F=0
\]
\end{maplelatex}

\begin{maplettyout}
\begin{center} {\bf Case   2}\end{center}
\end{maplettyout}

\begin{maplelatex}
\[
k=1, \,M=[], \,F= \left[ 
{
\begin{array}{r}
1
\end{array}}
 \right] 
\]
\end{maplelatex}

\begin{maplettyout}
using   case   1
\end{maplettyout}

\begin{maplettyout}
\begin{center}{\bf Case   3}\end{center}
\end{maplettyout}

\begin{maplelatex}
\[
k=2, \,M= \left[ 
{\begin{array}{r}
0
\end{array}}
 \right] , \,F= \left[ 
{\begin{array}{rr}
-1 & 1 \\
1 & 0
\end{array}}
 \right] 
\]
\end{maplelatex}

\begin{maplettyout}
covering   row   1
\end{maplettyout}

\begin{maplettyout}
using   cases   1   2
\end{maplettyout}

\begin{maplettyout}
\begin{center}{\bf Case   4}\end{center}
\end{maplettyout}

\begin{maplelatex}
\[
k=2, \,M=[], \,F= \left[ 
{\begin{array}{rr}
1 & 1 \\
1 & 0
\end{array}}
 \right] 
\]
\end{maplelatex}

\begin{maplettyout}
using   case   3
\end{maplettyout}

\begin{maplettyout}
\begin{center}{\bf Case   5}\end{center}
\end{maplettyout}

\begin{maplelatex}
\[
k=3, \,M= \left[ 
{\begin{array}{rc}
0 & 0 \\
0 & {X_{1}}
\end{array}}
 \right] , \,F= \left[ 
{\begin{array}{rrr}
0 & -1 & 1 \\
-1 & 1 & 0 \\
1 & 0 & 0
\end{array}}
 \right] 
\]
\end{maplelatex}

\begin{maplettyout}
covering   rows   1   2
\end{maplettyout}

\begin{maplettyout}
using   cases   1   2   3
\end{maplettyout}

\begin{maplettyout}
\begin{center}{\bf Case   6}\end{center}
\end{maplettyout}

\begin{maplelatex}
\[
k=3, \,M= \left[ 
{\begin{array}{r}
0 \\
0
\end{array}}
 \right] , \,F= \left[ 
{\begin{array}{rrr}
1 & 1 & 1 \\
-2 & 1 & 0 \\
1 & 0 & 0
\end{array}}
 \right] 
\]
\end{maplelatex}

\begin{maplettyout}
covering   column   1
\end{maplettyout}

\begin{maplettyout}
using   cases   3   5
\end{maplettyout}

\begin{maplettyout}
\begin{center}{\bf Case   7}\end{center}
\end{maplettyout}

\begin{maplelatex}
\[
k=3, \,M= \left[ 
{\begin{array}{cc}
{X_{1}} & 0 \\
0 & {X_{2}}
\end{array}}
 \right] , \,F= \left[ 
{\begin{array}{rrr}
1 & -2 & 1 \\
-2 & 2 & 0 \\
1 & 0 & 0
\end{array}}
 \right] 
\]
\end{maplelatex}

\begin{maplettyout}
covering   rows   1   2
\end{maplettyout}

\begin{maplettyout}
using   cases   1   3
\end{maplettyout}

\begin{maplettyout}
\begin{center}{\bf Case   8}\end{center}
\end{maplettyout}

\begin{maplelatex}
\[
k=3, \,M= \left[ 
{\begin{array}{r}
0
\end{array}}
 \right] , \,F= \left[ 
{\begin{array}{rrr}
-2 & 1 & 1 \\
1 & 1 & 0 \\
1 & 0 & 0
\end{array}}
 \right] 
\]
\end{maplelatex}

\begin{maplettyout}
covering   row   1
\end{maplettyout}

\begin{maplettyout}
using   cases   6   7
\end{maplettyout}

\begin{maplettyout}
\begin{center}{\bf Case   9}\end{center}
\end{maplettyout}

\begin{maplelatex}
\[
k=3, \,M=[], \,F= \left[ 
{\begin{array}{rrr}
1 & 1 & 1 \\
1 & 1 & 0 \\
1 & 0 & 0
\end{array}}
 \right] 
\]
\end{maplelatex}

\begin{maplettyout}
using   case   8
\end{maplettyout}

\begin{maplettyout}
\begin{center}{\bf Case   10}\end{center}
\end{maplettyout}

\begin{maplelatex}
\[
k=4, \,M= \left[ 
{\begin{array}{rcc}
0 & 0 & 0 \\
0 & 0 & {X_{1}} \\
0 & {X_{2}} & {X_{3}}
\end{array}}
 \right] , \,F= \left[ 
{\begin{array}{rrrr}
0 & 0 & -1 & 1 \\
0 & -1 & 1 & 0 \\
-1 & 1 & 0 & 0 \\
1 & 0 & 0 & 0
\end{array}}
 \right] 
\]
\end{maplelatex}

\begin{maplettyout}
covering   rows   1   2   3
\end{maplettyout}

\begin{maplettyout}
using   cases   1   2   3   5
\end{maplettyout}

\begin{maplettyout}
\begin{center}{\bf Case   11}\end{center}
\end{maplettyout}

\begin{maplelatex}
\[
k=4, \,M= \left[ 
{\begin{array}{rcc}
0 & 0 & 0 \\
0 & {X_{1}} & {X_{2}} \\
0 & {X_{3}} & {X_{4}}
\end{array}}
 \right] , \,F= \left[ 
{\begin{array}{rrrr}
0 & 1 & -2 & 1 \\
1 & 1 & 1 & 0 \\
-2 & 1 & 0 & 0 \\
1 & 0 & 0 & 0
\end{array}}
 \right] 
\]
\end{maplelatex}

\begin{maplettyout}
covering   row   1
\end{maplettyout}

\begin{maplettyout}
covering   column   1
\end{maplettyout}

\begin{maplettyout}
using   cases   3   5   10
\end{maplettyout}

\begin{maplettyout}
\begin{center}{\bf Case   12}\end{center}
\end{maplettyout}

\begin{maplelatex}
\[
k=4, \,M= \left[ 
{\begin{array}{rcc}
0 & {X_{1}} & 0 \\
0 & 0 & {X_{2}} \\
0 & {X_{3}} & {X_{4}}
\end{array}}
 \right] , \,F= \left[ 
{\begin{array}{rrrr}
0 & 1 & -2 & 1 \\
0 & -2 & 2 & 0 \\
-1 & 1 & 0 & 0 \\
1 & 0 & 0 & 0
\end{array}}
 \right] 
\]
\end{maplelatex}

\begin{maplettyout}
covering   rows   1   2   3
\end{maplettyout}

\begin{maplettyout}
using   cases   1   3   7
\end{maplettyout}

\begin{maplettyout}
\begin{center}{\bf Case   13}\end{center}
\end{maplettyout}

\begin{maplelatex}
\[
k=4, \,M= \left[ 
{\begin{array}{rc}
0 & 0 \\
0 & {X_{1}} \\
0 & {X_{2}}
\end{array}}
 \right] , \,F= \left[ 
{\begin{array}{rrrr}
0 & -2 & 1 & 1 \\
1 & 1 & 1 & 0 \\
-2 & 1 & 0 & 0 \\
1 & 0 & 0 & 0
\end{array}}
 \right] 
\]
\end{maplelatex}

\begin{maplettyout}
covering   columns   1   2
\end{maplettyout}

\begin{maplettyout}
using   cases   7   11   12
\end{maplettyout}

\begin{maplettyout}
\begin{center}{\bf Case   14}\end{center}
\end{maplettyout}

\begin{maplelatex}
\[
k=4, \,M= \left[ 
{\begin{array}{ccc}
0 & 0 & 0 \\
{X_{1}} & 0 & {X_{2}} \\
0 & {X_{3}} & {X_{4}}
\end{array}}
 \right] , \,F= \left[ 
{\begin{array}{rrrr}
0 & 0 & -1 & 1 \\
1 & -2 & 1 & 0 \\
-2 & 2 & 0 & 0 \\
1 & 0 & 0 & 0
\end{array}}
 \right] 
\]
\end{maplelatex}

\begin{maplettyout}
covering   rows   1   2   3
\end{maplettyout}

\begin{maplettyout}
using   cases   1   2   5
\end{maplettyout}

\begin{maplettyout}
\begin{center}{\bf Case   15}\end{center}
\end{maplettyout}

\begin{maplelatex}
\[
k=4, \,M= \left[ 
{\begin{array}{cc}
0 & 0 \\
{X_{1}} & 0 \\
0 & {X_{2}}
\end{array}}
 \right] , \,F= \left[ 
{\begin{array}{rrrr}
0 & 2 & 1 & 1 \\
1 & -4 & 1 & 0 \\
-2 & 2 & 0 & 0 \\
1 & 0 & 0 & 0
\end{array}}
 \right] 
\]
\end{maplelatex}

\begin{maplettyout}
covering   columns   1   2
\end{maplettyout}

\begin{maplettyout}
using   cases   3   10   14
\end{maplettyout}

\begin{maplettyout}
\begin{center}{\bf Case   16}\end{center}
\end{maplettyout}

\begin{maplelatex}
\[
k=4, \,M= \left[ 
{\begin{array}{ccc}
0 & {X_{1}} & 0 \\
{X_{2}} & 0 & {X_{3}} \\
0 & {X_{4}} & {X_{5}}
\end{array}}
 \right] , \,F= \left[ 
{\begin{array}{rrrr}
0 & 1 & -2 & 1 \\
1 & -3 & 2 & 0 \\
-2 & 2 & 0 & 0 \\
1 & 0 & 0 & 0
\end{array}}
 \right] 
\]
\end{maplelatex}

\begin{maplettyout}
covering   rows   1   2   3
\end{maplettyout}

\begin{maplettyout}
using   cases   1   3   5   7
\end{maplettyout}

\begin{maplettyout}
\begin{center}{\bf Case   17}\end{center}
\end{maplettyout}

\begin{maplelatex}
\[
k=4, \,M= \left[ 
{\begin{array}{rc}
0 & 0 \\
0 & {X_{1}}
\end{array}}
 \right] , \,F= \left[ 
{\begin{array}{rrrr}
-1 & -1 & 1 & 1 \\
-1 & 0 & 1 & 0 \\
1 & 1 & 0 & 0 \\
1 & 0 & 0 & 0
\end{array}}
 \right] 
\]
\end{maplelatex}

\begin{maplettyout}
covering   rows   1   2
\end{maplettyout}

\begin{maplettyout}
using   cases   13   15   16
\end{maplettyout}

\begin{maplettyout}
\begin{center}{\bf Case   18}\end{center}
\end{maplettyout}

\begin{maplelatex}
\[
k=4, \,M= \left[ 
{\begin{array}{cc}
{X_{1}} & 0 \\
0 & {X_{2}} \\
0 & {X_{3}}
\end{array}}
 \right] , \,F= \left[ 
{\begin{array}{rrrr}
-1 & -1 & 1 & 1 \\
3 & -1 & 1 & 0 \\
-3 & 2 & 0 & 0 \\
1 & 0 & 0 & 0
\end{array}}
 \right] 
\]
\end{maplelatex}

\begin{maplettyout}
covering   columns   1   2
\end{maplettyout}

\begin{maplettyout}
using   cases   7   11   16
\end{maplettyout}

\begin{maplettyout}
\begin{center}{\bf Case   19}\end{center}
\end{maplettyout}

\begin{maplelatex}
\[
k=4, \,M= \left[ 
{\begin{array}{r}
0 \\
0
\end{array}}
 \right] , \,F= \left[ 
{\begin{array}{rrrr}
3 & 1 & 1 & 1 \\
-5 & 1 & 1 & 0 \\
1 & 1 & 0 & 0 \\
1 & 0 & 0 & 0
\end{array}}
 \right] 
\]
\end{maplelatex}

\begin{maplettyout}
covering   column   1
\end{maplettyout}

\begin{maplettyout}
using   cases   17   18
\end{maplettyout}

\begin{maplettyout}
\begin{center}{\bf Case   20}\end{center}
\end{maplettyout}

\begin{maplelatex}
\[
k=4, \,M= \left[ 
{\begin{array}{ccc}
{X_{1}} & {X_{2}} & 0 \\
{X_{3}} & 0 & {X_{4}} \\
0 & {X_{5}} & {X_{6}}
\end{array}}
 \right] , \,F= \left[ 
{\begin{array}{rrrr}
-1 & 3 & -3 & 1 \\
3 & -6 & 3 & 0 \\
-3 & 3 & 0 & 0 \\
1 & 0 & 0 & 0
\end{array}}
 \right] 
\]
\end{maplelatex}

\begin{maplettyout}
covering   rows   1   2   3
\end{maplettyout}

\begin{maplettyout}
using   cases   1   7
\end{maplettyout}

\begin{maplettyout}
\begin{center}{\bf Case   21}\end{center}
\end{maplettyout}

\begin{maplelatex}
\[
k=4, \,M= \left[ 
{\begin{array}{cc}
{X_{1}} & 0 \\
0 & {X_{2}}
\end{array}}
 \right] , \,F= \left[ 
{\begin{array}{rrrr}
3 & -5 & 1 & 1 \\
-5 & 4 & 1 & 0 \\
1 & 1 & 0 & 0 \\
1 & 0 & 0 & 0
\end{array}}
 \right] 
\]
\end{maplelatex}

\begin{maplettyout}
covering   rows   1   2
\end{maplettyout}

\begin{maplettyout}
using   cases   18   20
\end{maplettyout}

\begin{maplettyout}
\begin{center}{\bf Case   22}\end{center}
\end{maplettyout}

\begin{maplelatex}
\[
k=4, \,M= \left[ 
{\begin{array}{r}
0
\end{array}}
 \right] , \,F= \left[ 
{\begin{array}{rrrr}
-3 & 1 & 1 & 1 \\
1 & 1 & 1 & 0 \\
1 & 1 & 0 & 0 \\
1 & 0 & 0 & 0
\end{array}}
 \right] 
\]
\end{maplelatex}

\begin{maplettyout}
covering   row   1
\end{maplettyout}

\begin{maplettyout}
using   cases   19   21
\end{maplettyout}

\begin{maplettyout}
\begin{center}{\bf Case   23}\end{center}
\end{maplettyout}

\begin{maplelatex}
\[
k=4, \,M=[], \,F= \left[ 
{\begin{array}{rrrr}
1 & 1 & 1 & 1 \\
1 & 1 & 1 & 0 \\
1 & 1 & 0 & 0 \\
1 & 0 & 0 & 0
\end{array}}
 \right] 
\]
\end{maplelatex}

\begin{maplettyout}
using   case   22
\end{maplettyout}

\begin{maplettyout}
\begin{center}{\bf Case   24}\end{center}
\end{maplettyout}

\begin{maplelatex}
\[
k=4, \,M= \left[ 
{\begin{array}{ccc}
0 & 0 & 0 \\
0 & 0 & {X_{1}} \\
0 & {X_{2}} & {X_{3}} \\
{\displaystyle \frac {{X_{4}}}{(m - 2)\,n}}  & {X_{5}} & {X_{6}}
\end{array}}
 \right] , \,F= \left[ 
{\begin{array}{rrrr}
0 & 0 & -1 & 1 \\
0 & -1 & 1 & 0 \\
0 & 1 & 0 & 0 \\
0 & 0 & 0 & 0
\end{array}}
 \right]  + {\displaystyle \frac {1}{(m - 2)\,(2\,m\,n - 5\,n - 1
)}} 
\]
\end{maplelatex}

\begin{maplettyout}
covering   rows   1   2   3
\end{maplettyout}

\begin{maplettyout}
using   cases   1   2   3   5
\end{maplettyout}

\begin{maplettyout}
\begin{center}{\bf Case   25}\end{center}
\end{maplettyout}

\begin{maplelatex}
\begin{eqnarray*}
\lefteqn{k=5, \,M= \left[ 
{\begin{array}{cccc}
0 & 0 & 0 & 0 \\
0 & 0 & 0 & {X_{1}} \\
{\displaystyle \frac {{X_{2}}}{(m - 2)\,(n - 1)}}  & 0 & {X_{3}}
 & {X_{4}} \\ [2ex]
0 & {X_{5}} & {X_{6}} & {X_{7}}
\end{array}}
 \right] , F= \left[ 
{\begin{array}{rrrrr}
0 & 0 & 0 & -1 & 1 \\
0 & 0 & -1 & 1 & 0 \\
0 & -1 & 1 & 0 & 0 \\
-1 & 1 & 0 & 0 & 0 \\
1 & 0 & 0 & 0 & 0
\end{array}}
 \right] } \\
 & & \mbox{} + {\displaystyle \frac {1}{2}} \,{\displaystyle 
\frac {2\,m - 5}{(m^{2} - 4\,m + 4)\,(2\,m\,n - 2\,m - 5\,n + 4)}
}  - {\displaystyle \frac {1}{2}} \,{\displaystyle \frac {1}{(m^{
2} - 4\,m + 4)\,n}} 
\end{eqnarray*}
\end{maplelatex}

\begin{maplettyout}
covering   rows   1   2   3   4
\end{maplettyout}

\begin{maplettyout}
using   cases   1   2   3   10   24
\end{maplettyout}

\begin{maplettyout}
\begin{center}{\bf Case   26}\end{center}
\end{maplettyout}

\begin{maplelatex}
\[
k=5, \,M= \left[ 
{\begin{array}{rccc}
0 & 0 & 0 & 0 \\
0 & 0 & 0 & {X_{1}} \\
0 & 0 & {X_{2}} & {X_{3}} \\
0 & {X_{4}} & {X_{5}} & {X_{6}}
\end{array}}
 \right] , \,F= \left[ 
{\begin{array}{rrrrr}
0 & 0 & 0 & -1 & 1 \\
0 & 0 & -1 & 1 & 0 \\
0 & -1 & 1 & 0 & 0 \\
-1 & 1 & 0 & 0 & 0 \\
1 & 0 & 0 & 0 & 0
\end{array}}
 \right] 
\]
\end{maplelatex}

\begin{maplettyout}
covering   rows   1   2   3   4
\end{maplettyout}

\begin{maplettyout}
using   cases   1   2   3   5   10
\end{maplettyout}

\begin{maplettyout}
\begin{center}{\bf Case   27}\end{center}
\end{maplettyout}

\begin{maplelatex}
\begin{eqnarray*}
\lefteqn{k=5, \,M= \left[ 
{\begin{array}{rccc}
0 & 0 & 0 & 0 \\
0 & 0 & 0 & {X_{1}} \\
0 & 0 & {X_{2}} & {X_{3}} \\
0 & {\displaystyle \frac {{X_{4}}}{(m - 2)\,(n - 1)}}  & {X_{5}}
 & {X_{6}}
\end{array}}
 \right] , F= \left[ 
{\begin{array}{rrrrr}
0 & 0 & 0 & -1 & 1 \\
0 & 0 & -1 & 1 & 0 \\
-1 & 0 & 1 & 0 & 0 \\
0 & 0 & 0 & 0 & 0 \\
1 & 0 & 0 & 0 & 0
\end{array}}
 \right] } \\
 & & \mbox{} + {\displaystyle \frac {1}{2}} \,{\displaystyle 
\frac {2\,m - 5}{(m^{2} - 4\,m + 4)\,(2\,m\,n - 2\,m - 5\,n + 4)}
}  - {\displaystyle \frac {1}{2}} \,{\displaystyle \frac {1}{(m^{
2} - 4\,m + 4)\,n}} 
\end{eqnarray*}
\end{maplelatex}

\begin{maplettyout}
covering   rows   1   2   3   4
\end{maplettyout}

\begin{maplettyout}
using   cases   1   2   3   5   24
\end{maplettyout}

\begin{maplettyout}
\begin{center}{\bf Case   28}\end{center}
\end{maplettyout}

\begin{maplelatex}
\[
k=3, \,M= \left[ 
{\begin{array}{cc}
0 & 0 \\
0 & {X_{1}} \\
{\displaystyle \frac {{X_{2}}}{(m - 1)\,n}}  & {X_{3}}
\end{array}}
 \right] , \,F= \left[ 
{\begin{array}{rrr}
0 & -1 & 1 \\
0 & 1 & 0 \\
0 & 0 & 0
\end{array}}
 \right]  + {\displaystyle \frac {1}{(m - 1)\,(2\,m\,n - 3\,n - 1
)}} 
\]
\end{maplelatex}

\begin{maplettyout}
covering   rows   1   2
\end{maplettyout}

\begin{maplettyout}
using   cases   1   2   3
\end{maplettyout}

\begin{maplettyout}
\begin{center}{\bf Case   29}\end{center}
\end{maplettyout}

\begin{maplelatex}
\begin{eqnarray*}
\lefteqn{k=4, \,M= \left[ 
{\begin{array}{ccc}
0 & 0 & 0 \\
{\displaystyle \frac {{X_{1}}}{(m - 1)\,(n - 1)}}  & 0 & {X_{2}}
 \\ [2ex]
0 & {X_{3}} & {X_{4}}
\end{array}}
 \right] , F=} \\
 & &  \left[ 
{\begin{array}{rrrr}
0 & 0 & -1 & 1 \\
0 & -1 & 1 & 0 \\
-1 & 1 & 0 & 0 \\
1 & 0 & 0 & 0
\end{array}}
 \right]  + {\displaystyle \frac {1}{2}} \,{\displaystyle \frac {
2\,m - 3}{(m^{2} - 2\,m + 1)\,(2\,m\,n - 2\,m - 3\,n + 2)}}  - 
{\displaystyle \frac {1}{2}} \,{\displaystyle \frac {1}{(m^{2} - 
2\,m + 1)\,n}} 
\end{eqnarray*}
\end{maplelatex}

\begin{maplettyout}
covering   rows   1   2   3
\end{maplettyout}

\begin{maplettyout}
using   cases   28   1   2   5
\end{maplettyout}

\begin{maplettyout}
\begin{center}{\bf Case   30}\end{center}
\end{maplettyout}

\begin{maplelatex}
\begin{eqnarray*}
\lefteqn{k=4, \,M= \left[ 
{\begin{array}{rcc}
0 & 0 & 0 \\
0 & 0 & {X_{1}} \\
0 & {\displaystyle \frac {{X_{2}}}{(m - 1)\,(n - 1)}}  & {X_{3}}
\end{array}}
 \right] , F=} \\
 & &  \left[ 
{\begin{array}{rrrr}
0 & 0 & -1 & 1 \\
-1 & 0 & 1 & 0 \\
0 & 0 & 0 & 0 \\
1 & 0 & 0 & 0
\end{array}}
 \right]  - {\displaystyle \frac {1}{2}} \,{\displaystyle \frac {
1}{(m^{2} - 2\,m + 1)\,n}}  + {\displaystyle \frac {1}{2}} \,
{\displaystyle \frac {2\,m - 3}{(m^{2} - 2\,m + 1)\,(2\,m\,n - 2
\,m - 3\,n + 2)}} 
\end{eqnarray*}
\end{maplelatex}

\begin{maplettyout}
covering   rows   1   2   3
\end{maplettyout}

\begin{maplettyout}
using   cases   28   1   2   3
\end{maplettyout}

\begin{maplettyout}
\begin{center}{\bf Case   31}\end{center}
\end{maplettyout}

\begin{maplelatex}
\begin{eqnarray*}
\lefteqn{k=5, \,M= \left[ 
{\begin{array}{cccc}
0 & 0 & 0 & 0 \\
0 & 0 & {X_{1}} & {X_{2}} \\
{\displaystyle \frac {{X_{3}}}{(m - 2)\,(n - 1)}}  & 0 & {X_{4}}
 & {X_{5}} \\ [2ex]
0 & {X_{6}} & {X_{7}} & {X_{8}}
\end{array}}
 \right] , F= \left[ 
{\begin{array}{rrrrr}
0 & 0 & 1 & -2 & 1 \\
0 & 1 & 1 & 1 & 0 \\
-1 & -2 & 1 & 0 & 0 \\
-1 & 1 & 0 & 0 & 0 \\
1 & 0 & 0 & 0 & 0
\end{array}}
 \right] } \\
 & & \mbox{} + {\displaystyle \frac {2\,m - 5}{(m^{2} - 4\,m + 4)
\,(m - 1)\,(2\,m\,n - 2\,m - 5\,n + 4)}}  + {\displaystyle 
\frac {m - 3}{(m^{2} - 4\,m + 4)\,n}} 
\end{eqnarray*}
\end{maplelatex}

\begin{maplettyout}
covering   row   1
\end{maplettyout}

\begin{maplettyout}
covering   columns   1   2
\end{maplettyout}

\begin{maplettyout}
using   cases   29   30   3   5   10   25   26   27
\end{maplettyout}

\begin{maplettyout}
\begin{center}{\bf Case   32}\end{center}
\end{maplettyout}

\begin{maplelatex}
\begin{eqnarray*}
\lefteqn{k=5, \,M= \left[ 
{\begin{array}{cccc}
0 & 0 & 0 & 0 \\
0 & {X_{1}} & 0 & {X_{2}} \\
{\displaystyle \frac {{X_{3}}}{(m - 2)\,(n - 1)}}  & 0 & {X_{4}}
 & {X_{5}} \\ [2ex]
0 & {X_{6}} & {X_{7}} & {X_{8}}
\end{array}}
 \right] , F= \left[ 
{\begin{array}{rrrrr}
0 & 0 & 0 & -1 & 1 \\
0 & 1 & -2 & 1 & 0 \\
0 & -2 & 2 & 0 & 0 \\
-1 & 1 & 0 & 0 & 0 \\
1 & 0 & 0 & 0 & 0
\end{array}}
 \right] } \\
 & & \mbox{} + {\displaystyle \frac {1}{2}} \,{\displaystyle 
\frac {2\,m - 5}{(m^{2} - 4\,m + 4)\,(2\,m\,n - 2\,m - 5\,n + 4)}
}  - {\displaystyle \frac {1}{2}} \,{\displaystyle \frac {1}{(m^{
2} - 4\,m + 4)\,n}} 
\end{eqnarray*}
\end{maplelatex}

\begin{maplettyout}
covering   rows   1   2   3   4
\end{maplettyout}

\begin{maplettyout}
using   cases   1   2   5   14   24
\end{maplettyout}

\begin{maplettyout}
\begin{center}{\bf Case   33}\end{center}
\end{maplettyout}

\begin{maplelatex}
\begin{eqnarray*}
\lefteqn{k=4, \,M= \left[ 
{\begin{array}{ccc}
0 & {X_{1}} & 0 \\
{\displaystyle \frac {{X_{2}}}{(m - 1)\,(n - 1)}}  & 0 & {X_{3}}
 \\ [2ex]
0 & {X_{4}} & {X_{5}}
\end{array}}
 \right] , F=} \\
 & &  \left[ 
{\begin{array}{rrrr}
0 & 1 & -2 & 1 \\
0 & -2 & 2 & 0 \\
-1 & 1 & 0 & 0 \\
1 & 0 & 0 & 0
\end{array}}
 \right]  + {\displaystyle \frac {1}{2}} \,{\displaystyle \frac {
2\,m - 3}{(m^{2} - 2\,m + 1)\,(2\,m\,n - 2\,m - 3\,n + 2)}}  - 
{\displaystyle \frac {1}{2}} \,{\displaystyle \frac {1}{(m^{2} - 
2\,m + 1)\,n}} 
\end{eqnarray*}
\end{maplelatex}

\begin{maplettyout}
covering   rows   1   2   3
\end{maplettyout}

\begin{maplettyout}
using   cases   28   1   3   7
\end{maplettyout}

\begin{maplettyout}
\begin{center}{\bf Case   34}\end{center}
\end{maplettyout}

\begin{maplelatex}
\begin{eqnarray*}
\lefteqn{k=5, \,M= \left[ 
{\begin{array}{cccc}
0 & 0 & 0 & 0 \\
{\displaystyle \frac {{X_{1}}}{(m - 2)\,(n - 1)}}  & 0 & {X_{2}}
 & {X_{3}} \\ [2ex]
{X_{4}} & 0 & {X_{5}} & {X_{6}} \\
0 & {X_{7}} & {X_{8}} & {X_{9}} \\
0 & {X_{10}} & {X_{11}} & {X_{12}}
\end{array}}
 \right] , } \\
 & & F= \left[ 
{\begin{array}{rrrrr}
0 & 0 & 1 & -2 & 1 \\
0 & 1 & 1 & 1 & 0 \\
0 & -2 & 1 & 0 & 0 \\
0 & 1 & 0 & 0 & 0 \\
0 & 0 & 0 & 0 & 0
\end{array}}
 \right]  + 2\,{\displaystyle \frac {1}{(2\,m\,n - 2\,m - 5\,n + 
4)\,(m - 1)\,(m - 2)}} 
\end{eqnarray*}
\end{maplelatex}

\begin{maplettyout}
covering   row   1
\end{maplettyout}

\begin{maplettyout}
covering   columns   1   2
\end{maplettyout}

\begin{maplettyout}
using   cases   28   3   5   10   24
\end{maplettyout}

\begin{maplettyout}
\begin{center}{\bf Case   35}\end{center}
\end{maplettyout}

\begin{maplelatex}
\[
k=5, \,M= \left[ 
{\begin{array}{rccc}
0 & 0 & 0 & 0 \\
0 & 0 & {X_{1}} & {X_{2}} \\
0 & {X_{3}} & {X_{4}} & {X_{5}} \\
0 & {X_{6}} & {X_{7}} & {X_{8}}
\end{array}}
 \right] , \,F= \left[ 
{\begin{array}{rrrrr}
0 & 0 & 1 & -2 & 1 \\
0 & -2 & 1 & 1 & 0 \\
1 & 1 & 1 & 0 & 0 \\
-2 & 1 & 0 & 0 & 0 \\
1 & 0 & 0 & 0 & 0
\end{array}}
 \right] 
\]
\end{maplelatex}

\begin{maplettyout}
covering   rows   1   2
\end{maplettyout}

\begin{maplettyout}
covering   column   1
\end{maplettyout}

\begin{maplettyout}
using   cases   31   32   33   34   7   14
\end{maplettyout}

\begin{maplettyout}
\begin{center}{\bf Case   36}\end{center}
\end{maplettyout}

\begin{maplelatex}
\[
k=5, \,M= \left[ 
{\begin{array}{rccc}
0 & 0 & {X_{1}} & 0 \\
0 & 0 & 0 & {X_{2}} \\
0 & 0 & {X_{3}} & {X_{4}} \\
0 & {X_{5}} & {X_{6}} & {X_{7}}
\end{array}}
 \right] , \,F= \left[ 
{\begin{array}{rrrrr}
0 & 0 & 1 & -2 & 1 \\
0 & 0 & -2 & 2 & 0 \\
0 & -1 & 1 & 0 & 0 \\
-1 & 1 & 0 & 0 & 0 \\
1 & 0 & 0 & 0 & 0
\end{array}}
 \right] 
\]
\end{maplelatex}

\begin{maplettyout}
covering   rows   1   2   3   4
\end{maplettyout}

\begin{maplettyout}
using   cases   1   3   7   12
\end{maplettyout}

\begin{maplettyout}
\begin{center}{\bf Case   37}\end{center}
\end{maplettyout}

\begin{maplelatex}
\[
k=5, \,M= \left[ 
{\begin{array}{rccc}
0 & 0 & {X_{1}} & 0 \\
0 & 0 & 0 & {X_{2}} \\
0 & {X_{3}} & {X_{4}} & {X_{5}} \\
0 & {X_{6}} & {X_{7}} & {X_{8}}
\end{array}}
 \right] , \,F= \left[ 
{\begin{array}{rrrrr}
0 & 0 & 1 & -2 & 1 \\
0 & 2 & -4 & 2 & 0 \\
1 & 1 & 1 & 0 & 0 \\
-2 & 1 & 0 & 0 & 0 \\
1 & 0 & 0 & 0 & 0
\end{array}}
 \right] 
\]
\end{maplelatex}

\begin{maplettyout}
covering   rows   1   2
\end{maplettyout}

\begin{maplettyout}
covering   column   1
\end{maplettyout}

\begin{maplettyout}
using   cases   36   3   5   10   12   26
\end{maplettyout}

\begin{maplettyout}
\begin{center}{\bf Case   38}\end{center}
\end{maplettyout}

\begin{maplelatex}
\[
k=5, \,M= \left[ 
{\begin{array}{rccc}
0 & 0 & {X_{1}} & 0 \\
0 & 0 & {X_{2}} & {X_{3}} \\
0 & {X_{4}} & 0 & {X_{5}} \\
0 & {X_{6}} & {X_{7}} & {X_{8}}
\end{array}}
 \right] , \,F= \left[ 
{\begin{array}{rrrrr}
0 & 0 & 1 & -2 & 1 \\
0 & 1 & -3 & 2 & 0 \\
0 & -2 & 2 & 0 & 0 \\
-1 & 1 & 0 & 0 & 0 \\
1 & 0 & 0 & 0 & 0
\end{array}}
 \right] 
\]
\end{maplelatex}

\begin{maplettyout}
covering   rows   1   2   3   4
\end{maplettyout}

\begin{maplettyout}
using   cases   1   3   5   7   16
\end{maplettyout}

\begin{maplettyout}
\begin{center}{\bf Case   39}\end{center}
\end{maplettyout}

\begin{maplelatex}
\[
k=5, \,M= \left[ 
{\begin{array}{rcc}
0 & 0 & 0 \\
0 & 0 & {X_{1}} \\
0 & {X_{2}} & {X_{3}} \\
0 & {X_{4}} & {X_{5}}
\end{array}}
 \right] , \,F= \left[ 
{\begin{array}{rrrrr}
0 & -1 & -1 & 1 & 1 \\
0 & -1 & 0 & 1 & 0 \\
1 & 1 & 1 & 0 & 0 \\
-2 & 1 & 0 & 0 & 0 \\
1 & 0 & 0 & 0 & 0
\end{array}}
 \right] 
\]
\end{maplelatex}

\begin{maplettyout}
covering   columns   1   2   3
\end{maplettyout}

\begin{maplettyout}
using   cases   35   37   38   16
\end{maplettyout}

\begin{maplettyout}
\begin{center}{\bf Case   40}\end{center}
\end{maplettyout}

\begin{maplelatex}
\[
k=5, \,M= \left[ 
{\begin{array}{cccc}
0 & 0 & 0 & 0 \\
0 & 0 & 0 & {X_{1}} \\
{X_{2}} & 0 & {X_{3}} & {X_{4}} \\
0 & {X_{5}} & {X_{6}} & {X_{7}}
\end{array}}
 \right] , \,F= \left[ 
{\begin{array}{rrrrr}
0 & 0 & 0 & -1 & 1 \\
0 & 0 & -1 & 1 & 0 \\
1 & -2 & 1 & 0 & 0 \\
-2 & 2 & 0 & 0 & 0 \\
1 & 0 & 0 & 0 & 0
\end{array}}
 \right] 
\]
\end{maplelatex}

\begin{maplettyout}
covering   rows   1   2   3   4
\end{maplettyout}

\begin{maplettyout}
using   cases   1   2   3   10
\end{maplettyout}

\begin{maplettyout}
\begin{center}{\bf Case   41}\end{center}
\end{maplettyout}

\begin{maplelatex}
\[
k=5, \,M= \left[ 
{\begin{array}{cccc}
0 & 0 & 0 & 0 \\
0 & 0 & {X_{1}} & {X_{2}} \\
{X_{3}} & 0 & {X_{4}} & {X_{5}} \\
0 & {X_{6}} & {X_{7}} & {X_{8}}
\end{array}}
 \right] , \,F= \left[ 
{\begin{array}{rrrrr}
0 & 0 & 1 & -2 & 1 \\
0 & 2 & 1 & 1 & 0 \\
1 & -4 & 1 & 0 & 0 \\
-2 & 2 & 0 & 0 & 0 \\
1 & 0 & 0 & 0 & 0
\end{array}}
 \right] 
\]
\end{maplelatex}

\begin{maplettyout}
covering   row   1
\end{maplettyout}

\begin{maplettyout}
covering   columns   1   2
\end{maplettyout}

\begin{maplettyout}
using   cases   40   3   5   10   14   26
\end{maplettyout}

\begin{maplettyout}
\begin{center}{\bf Case   42}\end{center}
\end{maplettyout}

\begin{maplelatex}
\[
k=3, \,M= \left[ 
{\begin{array}{cc}
0 & 0 \\
0 & {X_{1}} \\
{\displaystyle \frac {{X_{2}}}{(m - 1)\,(n + 1)}}  & {X_{3}}
\end{array}}
 \right] , \,F= \left[ 
{\begin{array}{rrr}
0 & -1 & 1 \\
0 & 1 & 0 \\
0 & 0 & 0
\end{array}}
 \right]  + {\displaystyle \frac {1}{(m - 1)\,(2\,m\,n - 3\,n - 2
 + m)}} 
\]
\end{maplelatex}

\begin{maplettyout}
covering   rows   1   2
\end{maplettyout}

\begin{maplettyout}
using   cases   1   2   3
\end{maplettyout}

\begin{maplettyout}
\begin{center}{\bf Case   43}\end{center}
\end{maplettyout}

\begin{maplelatex}
\begin{eqnarray*}
\lefteqn{k=4, \,M= \left[ 
{\begin{array}{ccc}
0 & {X_{1}} & 0 \\
{\displaystyle \frac {{X_{2}}}{(m - 1)\,n}}  & 0 & {X_{3}} \\
 [2ex]
0 & {X_{4}} & {X_{5}}
\end{array}}
 \right] , } \\
 & & F= \left[ 
{\begin{array}{rrrr}
0 & 1 & -2 & 1 \\
0 & -2 & 2 & 0 \\
-1 & 1 & 0 & 0 \\
1 & 0 & 0 & 0
\end{array}}
 \right]  + {\displaystyle \frac {2\,m - 3}{(m^{2} - 2\,m + 1)\,(
2\,m\,n - m - 3\,n + 1)}}  - {\displaystyle \frac {1}{(m^{2} - 2
\,m + 1)\,n}} 
\end{eqnarray*}
\end{maplelatex}

\begin{maplettyout}
covering   rows   1   2   3
\end{maplettyout}

\begin{maplettyout}
using   cases   42   1   3   7
\end{maplettyout}

\begin{maplettyout}
\begin{center}{\bf Case   44}\end{center}
\end{maplettyout}

\begin{maplelatex}
\begin{eqnarray*}
\lefteqn{k=5, \,M= \left[ 
{\begin{array}{cccc}
0 & 0 & {X_{1}} & 0 \\
0 & {\displaystyle \frac {{X_{2}}}{(m - 1)\,(n - 1)}}  & 0 & {X_{
3}} \\ [2ex]
{X_{4}} & 0 & {X_{5}} & {X_{6}} \\
0 & {X_{7}} & {X_{8}} & {X_{9}}
\end{array}}
 \right] , F= \left[ 
{\begin{array}{rrrrr}
0 & 0 & 1 & -2 & 1 \\
0 & 0 & -2 & 2 & 0 \\
1 & -2 & 1 & 0 & 0 \\
-2 & 2 & 0 & 0 & 0 \\
1 & 0 & 0 & 0 & 0
\end{array}}
 \right] } \\
 & & \mbox{} - {\displaystyle \frac {1}{(m^{2} - 2\,m + 1)\,(n - 
1)}}  + {\displaystyle \frac {2\,m - 3}{(m^{2} - 2\,m + 1)\,(2\,m
\,n - 3\,m - 3\,n + 4)}} 
\end{eqnarray*}
\end{maplelatex}

\begin{maplettyout}
covering   rows   1   2   3   4
\end{maplettyout}

\begin{maplettyout}
using   cases   42   43   1   3   12
\end{maplettyout}

\begin{maplettyout}
\begin{center}{\bf Case   45}\end{center}
\end{maplettyout}

\begin{maplelatex}
\begin{eqnarray*}
\lefteqn{k=5, \,M= \left[ 
{\begin{array}{rccc}
0 & 0 & {X_{1}} & 0 \\
0 & 0 & {X_{2}} & {X_{3}} \\
0 & {\displaystyle \frac {{X_{4}}}{(m - 1)\,(n - 1)}}  & 0 & {X_{
5}} \\ [2ex]
0 & {X_{6}} & {X_{7}} & {X_{8}}
\end{array}}
 \right] , F= \left[ 
{\begin{array}{rrrrr}
0 & 0 & 1 & -2 & 1 \\
0 & 0 & -2 & 2 & 0 \\
0 & -1 & 1 & 0 & 0 \\
-1 & 1 & 0 & 0 & 0 \\
1 & 0 & 0 & 0 & 0
\end{array}}
 \right]  - {\displaystyle \frac {1}{(m^{2} - 2\,m + 1)\,(n - 1)}
} } \\
 & & \mbox{} + {\displaystyle \frac {2\,m - 3}{(m^{2} - 2\,m + 1)
\,(2\,m\,n - 3\,m - 3\,n + 4)}} \mbox{\hspace{196pt}}
\end{eqnarray*}
\end{maplelatex}

\begin{maplettyout}
covering   rows   1   2   3   4
\end{maplettyout}

\begin{maplettyout}
using   cases   42   43   1   3   7
\end{maplettyout}

\begin{maplettyout}
\begin{center}{\bf Case   46}\end{center}
\end{maplettyout}

\begin{maplelatex}
\begin{eqnarray*}
\lefteqn{k=5, \,M= \left[ 
{\begin{array}{ccc}
0 & {\displaystyle \frac {{X_{1}}}{(m - 1)\,(n - 1)}}  & 0 \\
 [2ex]
0 & 0 & {X_{2}} \\
{X_{3}} & 0 & {X_{4}} \\
0 & {X_{5}} & {X_{6}}
\end{array}}
 \right] , F= \left[ 
{\begin{array}{rrrrr}
0 & -2 & -2 & 1 & 1 \\
0 & 2 & 1 & 1 & 0 \\
1 & -4 & 1 & 0 & 0 \\
-2 & 2 & 0 & 0 & 0 \\
1 & 0 & 0 & 0 & 0
\end{array}}
 \right] } \\
 & & \mbox{} + 2\,{\displaystyle \frac {2\,m - 3}{m\,(m^{2} - 2\,
m + 1)\,(2\,m\,n - 3\,m - 3\,n + 4)}}  + 2\,{\displaystyle 
\frac {m - 2}{(m^{2} - 2\,m + 1)\,(n - 1)}} 
\end{eqnarray*}
\end{maplelatex}

\begin{maplettyout}
covering   columns   1   2   3
\end{maplettyout}

\begin{maplettyout}
using   cases   36   41   44   45   7
\end{maplettyout}

\begin{maplettyout}
\begin{center}{\bf Case   47}\end{center}
\end{maplettyout}

\begin{maplelatex}
\[
k=5, \,M= \left[ 
{\begin{array}{cccc}
0 & 0 & 0 & 0 \\
0 & 0 & {X_{1}} & {X_{2}} \\
{X_{3}} & {X_{4}} & 0 & {X_{5}} \\
0 & {X_{6}} & {X_{7}} & {X_{8}}
\end{array}}
 \right] , \,F= \left[ 
{\begin{array}{rrrrr}
0 & 0 & 0 & -1 & 1 \\
0 & 1 & -2 & 1 & 0 \\
1 & -3 & 2 & 0 & 0 \\
-2 & 2 & 0 & 0 & 0 \\
1 & 0 & 0 & 0 & 0
\end{array}}
 \right] 
\]
\end{maplelatex}

\begin{maplettyout}
covering   rows   1   2   3   4
\end{maplettyout}

\begin{maplettyout}
using   cases   1   2   5   10   14
\end{maplettyout}

\begin{maplettyout}
\begin{center}{\bf Case   48}\end{center}
\end{maplettyout}

\begin{maplelatex}
\begin{eqnarray*}
\lefteqn{k=4, \,M= \left[ 
{\begin{array}{ccc}
0 & 0 & 0 \\
{X_{1}} & 0 & {X_{2}} \\
0 & {\displaystyle \frac {{X_{3}}}{(m - 1)\,n}}  & {X_{4}}
\end{array}}
 \right] , } \\
 & & F= \left[ 
{\begin{array}{rrrr}
0 & 0 & -1 & 1 \\
0 & -1 & 1 & 0 \\
-1 & 1 & 0 & 0 \\
1 & 0 & 0 & 0
\end{array}}
 \right]  + {\displaystyle \frac {2\,m - 3}{(m^{2} - 2\,m + 1)\,(
2\,m\,n - m - 3\,n + 1)}}  - {\displaystyle \frac {1}{(m^{2} - 2
\,m + 1)\,n}} 
\end{eqnarray*}
\end{maplelatex}

\begin{maplettyout}
covering   rows   1   2   3
\end{maplettyout}

\begin{maplettyout}
using   cases   42   1   2   5
\end{maplettyout}

\begin{maplettyout}
\begin{center}{\bf Case   49}\end{center}
\end{maplettyout}

\begin{maplelatex}
\begin{eqnarray*}
\lefteqn{k=5, \,M= \left[ 
{\begin{array}{cccc}
0 & 0 & 0 & 0 \\
0 & 0 & {\displaystyle \frac {{X_{1}}}{(m - 1)\,(n - 1)}}  & {X_{
2}} \\ [2ex]
{X_{3}} & {X_{4}} & 0 & {X_{5}} \\
0 & {X_{6}} & {X_{7}} & {X_{8}}
\end{array}}
 \right] , F= \left[ 
{\begin{array}{rrrrr}
0 & 0 & 0 & -1 & 1 \\
0 & 0 & -1 & 1 & 0 \\
1 & -2 & 1 & 0 & 0 \\
-2 & 2 & 0 & 0 & 0 \\
1 & 0 & 0 & 0 & 0
\end{array}}
 \right] } \\
 & & \mbox{} + {\displaystyle \frac {2\,m - 3}{(m^{2} - 2\,m + 1)
\,(2\,m\,n - 3\,m - 3\,n + 4)}}  - {\displaystyle \frac {1}{(m^{2
} - 2\,m + 1)\,(n - 1)}} 
\end{eqnarray*}
\end{maplelatex}

\begin{maplettyout}
covering   rows   1   2   3   4
\end{maplettyout}

\begin{maplettyout}
using   cases   42   48   1   2   10
\end{maplettyout}

\begin{maplettyout}
\begin{center}{\bf Case   50}\end{center}
\end{maplettyout}

\begin{maplelatex}
\begin{eqnarray*}
\lefteqn{k=4, \,M= \left[ 
{\begin{array}{rcc}
0 & 0 & 0 \\
0 & 0 & {X_{1}} \\
0 & {\displaystyle \frac {{X_{2}}}{(m - 1)\,n}}  & {X_{3}}
\end{array}}
 \right] , } \\
 & & F= \left[ 
{\begin{array}{rrrr}
0 & 0 & -1 & 1 \\
-1 & 0 & 1 & 0 \\
0 & 0 & 0 & 0 \\
1 & 0 & 0 & 0
\end{array}}
 \right]  - {\displaystyle \frac {1}{(m^{2} - 2\,m + 1)\,n}}  + 
{\displaystyle \frac {2\,m - 3}{(m^{2} - 2\,m + 1)\,(2\,m\,n - m
 - 3\,n + 1)}} 
\end{eqnarray*}
\end{maplelatex}

\begin{maplettyout}
covering   rows   1   2   3
\end{maplettyout}

\begin{maplettyout}
using   cases   42   1   2   3
\end{maplettyout}

\begin{maplettyout}
\begin{center}{\bf Case   51}\end{center}
\end{maplettyout}

\begin{maplelatex}
\begin{eqnarray*}
\lefteqn{k=5, \,M= \left[ 
{\begin{array}{rccc}
0 & 0 & 0 & 0 \\
0 & 0 & 0 & {X_{1}} \\
0 & 0 & {\displaystyle \frac {{X_{2}}}{(m - 1)\,(n - 1)}}  & {X_{
3}} \\ [2ex]
0 & {X_{4}} & {X_{5}} & {X_{6}}
\end{array}}
 \right] , F= \left[ 
{\begin{array}{rrrrr}
0 & 0 & 0 & -1 & 1 \\
0 & -1 & 0 & 1 & 0 \\
0 & 0 & 0 & 0 & 0 \\
-1 & 1 & 0 & 0 & 0 \\
1 & 0 & 0 & 0 & 0
\end{array}}
 \right] } \\
 & & \mbox{} + {\displaystyle \frac {2\,m - 3}{(m^{2} - 2\,m + 1)
\,(2\,m\,n - 3\,m - 3\,n + 4)}}  - {\displaystyle \frac {1}{(m^{2
} - 2\,m + 1)\,(n - 1)}} 
\end{eqnarray*}
\end{maplelatex}

\begin{maplettyout}
covering   rows   1   2   3   4
\end{maplettyout}

\begin{maplettyout}
using   cases   42   50   1   2   3
\end{maplettyout}

\begin{maplettyout}
\begin{center}{\bf Case   52}\end{center}
\end{maplettyout}

\begin{maplelatex}
\begin{eqnarray*}
\lefteqn{k=5, \,M= \left[ 
{\begin{array}{ccc}
{\displaystyle \frac {{X_{1}}}{(m - 1)\,(n - 1)}}  & 0 & 0 \\
 [2ex]
0 & 0 & {X_{2}} \\
{X_{3}} & {X_{4}} & 0 \\
0 & {X_{5}} & {X_{6}}
\end{array}}
 \right] , F=} \\
 & &  \left[ 
{\begin{array}{rrrrr}
0 & -1 & 2 & 1 & 1 \\
0 & 0 & -4 & 1 & 0 \\
1 & -3 & 2 & 0 & 0 \\
-2 & 2 & 0 & 0 & 0 \\
1 & 0 & 0 & 0 & 0
\end{array}}
 \right]  + 2\,{\displaystyle \frac {m - 2}{(m - 1)^{2}\,(n - 1)}
}  + 2\,{\displaystyle \frac {2\,m - 3}{m\,(m - 1)^{2}\,(2\,m\,n
 - 3\,m - 3\,n + 4)}} 
\end{eqnarray*}
\end{maplelatex}

\begin{maplettyout}
covering   columns   1   2   3
\end{maplettyout}

\begin{maplettyout}
using   cases   47   49   51   3   26
\end{maplettyout}

\begin{maplettyout}
\begin{center}{\bf Case   53}\end{center}
\end{maplettyout}

\begin{maplelatex}
\[
k=5, \,M= \left[ 
{\begin{array}{cccc}
0 & 0 & {X_{1}} & 0 \\
0 & 0 & {X_{2}} & {X_{3}} \\
{X_{4}} & {X_{5}} & 0 & {X_{6}} \\
0 & {X_{7}} & {X_{8}} & {X_{9}}
\end{array}}
 \right] , \,F= \left[ 
{\begin{array}{rrrrr}
0 & 0 & 1 & -2 & 1 \\
0 & 1 & -3 & 2 & 0 \\
1 & -3 & 2 & 0 & 0 \\
-2 & 2 & 0 & 0 & 0 \\
1 & 0 & 0 & 0 & 0
\end{array}}
 \right] 
\]
\end{maplelatex}

\begin{maplettyout}
covering   rows   1   2   3   4
\end{maplettyout}

\begin{maplettyout}
using   cases   1   3   7   10   16
\end{maplettyout}

\begin{maplettyout}
\begin{center}{\bf Case   54}\end{center}
\end{maplettyout}

\begin{maplelatex}
\begin{eqnarray*}
\lefteqn{k=5, \,M= \left[ 
{\begin{array}{ccc}
{\displaystyle \frac {{X_{1}}}{(m - 1)\,(n - 1)}}  & 0 & 0 \\
 [2ex]
0 & 0 & {X_{2}} \\
0 & {X_{3}} & {X_{4}}
\end{array}}
 \right] , F= \left[ 
{\begin{array}{rrrrr}
0 & -4 & -1 & 1 & 1 \\
-4 & 3 & 0 & 1 & 0 \\
-1 & 0 & 1 & 0 & 0 \\
1 & 1 & 0 & 0 & 0 \\
1 & 0 & 0 & 0 & 0
\end{array}}
 \right]  + 9\,{\displaystyle \frac {1}{(3\,m - 4)\,n}} } \\
 & & \mbox{} - 4\,{\displaystyle \frac {1}{(m^{2} - 2\,m + 1)\,(n
 - 1)}}  \\
 & & \mbox{} + 4\,{\displaystyle \frac {(2\,m - 3)^{2}}{(m^{2} - 
2\,m + 1)\,m\,(3\,m - 4)\,(2\,m\,n - 3\,m - 3\,n + 4)}} 
\mbox{\hspace{63pt}}
\end{eqnarray*}
\end{maplelatex}

\begin{maplettyout}
covering   rows   1   2   3
\end{maplettyout}

\begin{maplettyout}
using   cases   39   46   52   53
\end{maplettyout}

\begin{maplettyout}
\begin{center}{\bf Case   55}\end{center}
\end{maplettyout}

\begin{maplelatex}
\begin{eqnarray*}
\lefteqn{k=5, \,M= \left[ 
{\begin{array}{cccc}
0 & 0 & {X_{1}} & 0 \\
0 & 0 & {X_{2}} & {X_{3}} \\
{\displaystyle \frac {{X_{4}}}{(m - 2)\,(n - 1)}}  & {X_{5}} & 0
 & {X_{6}} \\ [2ex]
0 & {X_{7}} & {X_{8}} & {X_{9}}
\end{array}}
 \right] , F= \left[ 
{\begin{array}{rrrrr}
0 & 0 & 1 & -2 & 1 \\
0 & 1 & -3 & 2 & 0 \\
0 & -2 & 2 & 0 & 0 \\
-1 & 1 & 0 & 0 & 0 \\
1 & 0 & 0 & 0 & 0
\end{array}}
 \right] } \\
 & & \mbox{} + {\displaystyle \frac {1}{2}} \,{\displaystyle 
\frac {2\,m - 5}{(m^{2} - 4\,m + 4)\,(2\,m\,n - 2\,m - 5\,n + 4)}
}  - {\displaystyle \frac {1}{2}} \,{\displaystyle \frac {1}{(m^{
2} - 4\,m + 4)\,n}} 
\end{eqnarray*}
\end{maplelatex}

\begin{maplettyout}
covering   rows   1   2   3   4
\end{maplettyout}

\begin{maplettyout}
using   cases   1   3   7   16   24
\end{maplettyout}

\begin{maplettyout}
\begin{center}{\bf Case   56}\end{center}
\end{maplettyout}

\begin{maplelatex}
\[
k=5, \,M= \left[ 
{\begin{array}{rccc}
0 & {X_{1}} & 0 & 0 \\
0 & 0 & {X_{2}} & {X_{3}} \\
0 & {X_{4}} & {X_{5}} & {X_{6}} \\
0 & {X_{7}} & {X_{8}} & {X_{9}}
\end{array}}
 \right] , \,F= \left[ 
{\begin{array}{rrrrr}
0 & -1 & 3 & -3 & 1 \\
0 & -1 & -1 & 2 & 0 \\
1 & 1 & 1 & 0 & 0 \\
-2 & 1 & 0 & 0 & 0 \\
1 & 0 & 0 & 0 & 0
\end{array}}
 \right] 
\]
\end{maplelatex}

\begin{maplettyout}
covering   rows   1   2
\end{maplettyout}

\begin{maplettyout}
covering   column   1
\end{maplettyout}

\begin{maplettyout}
using   cases   31   33   34   55   7   16
\end{maplettyout}

\begin{maplettyout}
\begin{center}{\bf Case   57}\end{center}
\end{maplettyout}

\begin{maplelatex}
\begin{eqnarray*}
\lefteqn{k=5, \,M= \left[ 
{\begin{array}{cccc}
{\displaystyle \frac {{X_{1}}}{(m - 1)\,(n - 1)}}  & 0 & {X_{2}}
 & 0 \\ [2ex]
0 & {X_{3}} & 0 & {X_{4}} \\
{X_{5}} & 0 & {X_{6}} & {X_{7}} \\
0 & {X_{8}} & {X_{9}} & {X_{10}}
\end{array}}
 \right] , F= \left[ 
{\begin{array}{rrrrr}
0 & 0 & 1 & -2 & 1 \\
0 & 1 & -3 & 2 & 0 \\
1 & -3 & 2 & 0 & 0 \\
-2 & 2 & 0 & 0 & 0 \\
1 & 0 & 0 & 0 & 0
\end{array}}
 \right] } \\
 & & \mbox{} - {\displaystyle \frac {1}{(m^{2} - 2\,m + 1)\,(n - 
1)}}  + {\displaystyle \frac {2\,m - 3}{(m^{2} - 2\,m + 1)\,(2\,m
\,n - 3\,m - 3\,n + 4)}} 
\end{eqnarray*}
\end{maplelatex}

\begin{maplettyout}
covering   rows   1   2   3   4
\end{maplettyout}

\begin{maplettyout}
using   cases   42   43   1   5   16
\end{maplettyout}

\begin{maplettyout}
\begin{center}{\bf Case   58}\end{center}
\end{maplettyout}

\begin{maplelatex}
\begin{eqnarray*}
\lefteqn{k=5, \,M= \left[ 
{\begin{array}{cccc}
{\displaystyle \frac {{X_{1}}}{(m - 1)\,(n - 1)}}  & {X_{2}} & 0
 & 0 \\ [2ex]
0 & 0 & {X_{3}} & {X_{4}} \\
0 & {X_{5}} & {X_{6}} & {X_{7}}
\end{array}}
 \right] , F= \left[ 
{\begin{array}{rrrrr}
0 & -1 & 3 & -3 & 1 \\
-3 & 2 & -1 & 2 & 0 \\
-1 & 0 & 1 & 0 & 0 \\
1 & 1 & 0 & 0 & 0 \\
1 & 0 & 0 & 0 & 0
\end{array}}
 \right] } \\
 & & \mbox{} + 2\,{\displaystyle \frac {(2\,m - 3)^{2}}{(3\,m - 4
)\,(m - 1)^{2}\,(2\,m\,n - 3\,m - 3\,n + 4)}}  - 2\,
{\displaystyle \frac {m}{(m - 1)^{2}\,(n - 1)}}  \\
 & & \mbox{} + 6\,{\displaystyle \frac {1}{(3\,m - 4)\,n}} 
\end{eqnarray*}
\end{maplelatex}

\begin{maplettyout}
covering   rows   1   2   3
\end{maplettyout}

\begin{maplettyout}
using   cases   41   56   57   16
\end{maplettyout}

\begin{maplettyout}
\begin{center}{\bf Case   59}\end{center}
\end{maplettyout}

\begin{maplelatex}
\begin{eqnarray*}
\lefteqn{k=4, \,M= \left[ 
{\begin{array}{ccc}
0 & 0 & 0 \\
0 & 0 & {X_{1}} \\
0 & {X_{2}} & {X_{3}} \\
{X_{4}} + {\displaystyle \frac {{X_{5}}}{(m - 2)\,n}}  & {X_{6}}
 & {X_{7}}
\end{array}}
 \right] , } \\
 & & F= \left[ 
{\begin{array}{rrrr}
0 & 0 & -1 & 1 \\
0 & -1 & 1 & 0 \\
0 & 1 & 0 & 0 \\
0 & 0 & 0 & 0
\end{array}}
 \right]  + {\displaystyle \frac {2\,m - 5}{(6 - 5\,m + m^{2})\,(
2\,m\,n - 5\,n - 1)}} 
\end{eqnarray*}
\end{maplelatex}

\begin{maplettyout}
covering   rows   1   2   3
\end{maplettyout}

\begin{maplettyout}
using   cases   1   2   3   5   24
\end{maplettyout}

\begin{maplettyout}
\begin{center}{\bf Case   60}\end{center}
\end{maplettyout}

\begin{maplelatex}
\begin{eqnarray*}
\lefteqn{k=5, \,M= \left[ 
{\begin{array}{cccc}
0 & 0 & 0 & 0 \\
0 & 0 & 0 & {X_{1}} \\
{X_{2}} + {\displaystyle \frac {{X_{3}}}{(m - 2)\,(n - 1)}}  & 0
 & {X_{4}} & {X_{5}} \\ [2ex]
0 & {X_{6}} & {X_{7}} & {X_{8}}
\end{array}}
 \right] , F= \left[ 
{\begin{array}{rrrrr}
0 & 0 & 0 & -1 & 1 \\
0 & 0 & -1 & 1 & 0 \\
0 & -1 & 1 & 0 & 0 \\
0 & 1 & 0 & 0 & 0 \\
1 & 0 & 0 & 0 & 0
\end{array}}
 \right] } \\
 & & \mbox{} + {\displaystyle \frac {1}{2}} \,{\displaystyle 
\frac {(2\,m - 5)^{2}}{(16\,m - 7\,m^{2} + m^{3} - 12)\,(2\,m\,n
 - 2\,m - 5\,n + 4)}}  + {\displaystyle \frac {1}{2}} \,
{\displaystyle \frac {m - 3}{(m^{2} - 4\,m + 4)\,n}}  \\
 & & \mbox{} - {\displaystyle \frac {3}{2}} \,{\displaystyle 
\frac {1}{(m - 3)\,n}} 
\end{eqnarray*}
\end{maplelatex}

\begin{maplettyout}
covering   rows   1   2   3   4
\end{maplettyout}

\begin{maplettyout}
using   cases   59   1   2   3   10
\end{maplettyout}

\begin{maplettyout}
\begin{center}{\bf Case   61}\end{center}
\end{maplettyout}

\begin{maplelatex}
\begin{eqnarray*}
\lefteqn{k=5, \,M= \left[ 
{\begin{array}{rccc}
0 & 0 & 0 & 0 \\
0 & 0 & 0 & {X_{1}} \\
0 & 0 & {X_{2}} & {X_{3}} \\
0 & {X_{4}} + {\displaystyle \frac {{X_{5}}}{(m - 2)\,(n - 1)}} 
 & {X_{6}} & {X_{7}}
\end{array}}
 \right] , F= \left[ 
{\begin{array}{rrrrr}
0 & 0 & 0 & -1 & 1 \\
0 & 0 & -1 & 1 & 0 \\
-1 & 0 & 1 & 0 & 0 \\
0 & 0 & 0 & 0 & 0 \\
1 & 0 & 0 & 0 & 0
\end{array}}
 \right]  - {\displaystyle \frac {1}{2}} \,{\displaystyle \frac {
1}{(m - 3)\,n}} } \\
 & & \mbox{} + {\displaystyle \frac {1}{2}} \,{\displaystyle 
\frac {m - 3}{(m^{2} - 4\,m + 4)\,n}}  + {\displaystyle \frac {1
}{2}} \,{\displaystyle \frac {(2\,m - 5)^{2}}{(m^{2} - 4\,m + 4)
\,(m - 3)\,(2\,m\,n - 2\,m - 5\,n + 4)}} \mbox{\hspace{14pt}}
\end{eqnarray*}
\end{maplelatex}

\begin{maplettyout}
covering   rows   1   2   3   4
\end{maplettyout}

\begin{maplettyout}
using   cases   59   1   2   3   5
\end{maplettyout}

\begin{maplettyout}
\begin{center}{\bf Case   62}\end{center}
\end{maplettyout}

\begin{maplelatex}
\[
k=3, \,M= \left[ 
{\begin{array}{cc}
0 & 0 \\
0 & {X_{1}} \\
{X_{2}} + {\displaystyle \frac {{X_{3}}}{(m - 1)\,n}}  & {X_{4}}
\end{array}}
 \right] , \,F= \left[ 
{\begin{array}{rrr}
0 & -1 & 1 \\
0 & 1 & 0 \\
0 & 0 & 0
\end{array}}
 \right]  + {\displaystyle \frac {2\,m - 3}{(2\,m\,n - 3\,n - 1)
\,(m - 1)\,(m - 2)}} 
\]
\end{maplelatex}

\begin{maplettyout}
covering   rows   1   2
\end{maplettyout}

\begin{maplettyout}
using   cases   28   1   2   3
\end{maplettyout}

\begin{maplettyout}
\begin{center}{\bf Case   63}\end{center}
\end{maplettyout}

\begin{maplelatex}
\begin{eqnarray*}
\lefteqn{k=4, \,M= \left[ 
{\begin{array}{ccc}
0 & 0 & 0 \\
{X_{1}} + {\displaystyle \frac {{X_{2}}}{(m - 1)\,(n - 1)}}  & 0
 & {X_{3}} \\ [2ex]
0 & {X_{4}} & {X_{5}}
\end{array}}
 \right] , F= \left[ 
{\begin{array}{rrrr}
0 & 0 & -1 & 1 \\
0 & -1 & 1 & 0 \\
0 & 1 & 0 & 0 \\
1 & 0 & 0 & 0
\end{array}}
 \right]  - {\displaystyle \frac {3}{2}} \,{\displaystyle \frac {
1}{(m - 2)\,n}} } \\
 & & \mbox{} + {\displaystyle \frac {1}{2}} \,{\displaystyle 
\frac {(2\,m - 3)^{2}}{(m^{2} - 2\,m + 1)\,(m - 2)\,(2\,m\,n - 2
\,m - 3\,n + 2)}}  + {\displaystyle \frac {1}{2}} \,
{\displaystyle \frac {m - 2}{(m^{2} - 2\,m + 1)\,n}} 
\end{eqnarray*}
\end{maplelatex}

\begin{maplettyout}
covering   rows   1   2   3
\end{maplettyout}

\begin{maplettyout}
using   cases   62   1   2   5
\end{maplettyout}

\begin{maplettyout}
\begin{center}{\bf Case   64}\end{center}
\end{maplettyout}

\begin{maplelatex}
\begin{eqnarray*}
\lefteqn{k=4, \,M= \left[ 
{\begin{array}{rcc}
0 & 0 & 0 \\
0 & 0 & {X_{1}} \\
0 & {X_{2}} + {\displaystyle \frac {{X_{3}}}{(m - 1)\,(n - 1)}} 
 & {X_{4}}
\end{array}}
 \right] , F= \left[ 
{\begin{array}{rrrr}
0 & 0 & -1 & 1 \\
-1 & 0 & 1 & 0 \\
0 & 0 & 0 & 0 \\
1 & 0 & 0 & 0
\end{array}}
 \right]  - {\displaystyle \frac {1}{2}} \,{\displaystyle \frac {
1}{(m - 2)\,n}} } \\
 & & \mbox{} + {\displaystyle \frac {1}{2}} \,{\displaystyle 
\frac {m - 2}{(m^{2} - 2\,m + 1)\,n}}  + {\displaystyle \frac {1
}{2}} \,{\displaystyle \frac {(2\,m - 3)^{2}}{(m - 2)\,(m^{2} - 2
\,m + 1)\,(2\,m\,n - 2\,m - 3\,n + 2)}} 
\end{eqnarray*}
\end{maplelatex}

\begin{maplettyout}
covering   rows   1   2   3
\end{maplettyout}

\begin{maplettyout}
using   cases   62   1   2   3
\end{maplettyout}

\begin{maplettyout}
\begin{center}{\bf Case   65}\end{center}
\end{maplettyout}

\begin{maplelatex}
\begin{eqnarray*}
\lefteqn{k=5, \,M= \left[ 
{\begin{array}{cccc}
0 & 0 & 0 & 0 \\
0 & 0 & {X_{1}} & {X_{2}} \\
{X_{3}} + {\displaystyle \frac {{X_{4}}}{(m - 2)\,(n - 1)}}  & 0
 & {X_{5}} & {X_{6}} \\ [2ex]
0 & {X_{7}} & {X_{8}} & {X_{9}}
\end{array}}
 \right] , F= \left[ 
{\begin{array}{rrrrr}
0 & 0 & 1 & -2 & 1 \\
0 & 1 & 1 & 1 & 0 \\
-1 & -2 & 1 & 0 & 0 \\
0 & 1 & 0 & 0 & 0 \\
1 & 0 & 0 & 0 & 0
\end{array}}
 \right]  - {\displaystyle \frac {3}{2}} \,{\displaystyle \frac {
1}{(m - 3)\,n}} } \\
 & & \mbox{} - {\displaystyle \frac {1}{2}} \,{\displaystyle 
\frac {1}{(m - 1)\,n}}  + {\displaystyle \frac {2\,m - 5}{(m^{2}
 - 4\,m + 4)\,n}}  \\
 & & \mbox{} + {\displaystyle \frac {(2\,m - 5)^{2}}{(m^{2} - 4\,
m + 4)\,(m - 1)\,(m - 3)\,(2\,m\,n - 2\,m - 5\,n + 4)}} 
\mbox{\hspace{93pt}}
\end{eqnarray*}
\end{maplelatex}

\begin{maplettyout}
covering   row   1
\end{maplettyout}

\begin{maplettyout}
covering   columns   1   2
\end{maplettyout}

\begin{maplettyout}
using   cases   60   61   63   64   3   5   10   26
\end{maplettyout}

\begin{maplettyout}
\begin{center}{\bf Case   66}\end{center}
\end{maplettyout}

\begin{maplelatex}
\begin{eqnarray*}
\lefteqn{k=5, \,M= \left[ 
{\begin{array}{cccc}
0 & 0 & 0 & 0 \\
0 & {X_{1}} & 0 & {X_{2}} \\
{X_{3}} + {\displaystyle \frac {{X_{4}}}{(m - 2)\,(n - 1)}}  & 0
 & {X_{5}} & {X_{6}} \\ [2ex]
0 & {X_{7}} & {X_{8}} & {X_{9}}
\end{array}}
 \right] , F= \left[ 
{\begin{array}{rrrrr}
0 & 0 & 0 & -1 & 1 \\
0 & 1 & -2 & 1 & 0 \\
0 & -2 & 2 & 0 & 0 \\
0 & 1 & 0 & 0 & 0 \\
1 & 0 & 0 & 0 & 0
\end{array}}
 \right] } \\
 & & \mbox{} + {\displaystyle \frac {1}{2}} \,{\displaystyle 
\frac {(2\,m - 5)^{2}}{(16\,m - 7\,m^{2} + m^{3} - 12)\,(2\,m\,n
 - 2\,m - 5\,n + 4)}}  + {\displaystyle \frac {1}{2}} \,
{\displaystyle \frac {m - 3}{(m^{2} - 4\,m + 4)\,n}}  \\
 & & \mbox{} - {\displaystyle \frac {3}{2}} \,{\displaystyle 
\frac {1}{(m - 3)\,n}} 
\end{eqnarray*}
\end{maplelatex}

\begin{maplettyout}
covering   rows   1   2   3   4
\end{maplettyout}

\begin{maplettyout}
using   cases   59   1   2   5   14
\end{maplettyout}

\begin{maplettyout}
\begin{center}{\bf Case   67}\end{center}
\end{maplettyout}

\begin{maplelatex}
\begin{eqnarray*}
\lefteqn{k=4, \,M= \left[ 
{\begin{array}{ccc}
0 & {X_{1}} & 0 \\
{X_{2}} + {\displaystyle \frac {{X_{3}}}{(m - 1)\,(n - 1)}}  & 0
 & {X_{4}} \\ [2ex]
0 & {X_{5}} & {X_{6}}
\end{array}}
 \right] , F= \left[ 
{\begin{array}{rrrr}
0 & 1 & -2 & 1 \\
0 & -2 & 2 & 0 \\
0 & 1 & 0 & 0 \\
1 & 0 & 0 & 0
\end{array}}
 \right] } \\
 & & \mbox{} + {\displaystyle \frac {1}{2}} \,{\displaystyle 
\frac {(2\,m - 3)^{2}}{(m - 2)\,(m^{2} - 2\,m + 1)\,(2\,m\,n - 2
\,m - 3\,n + 2)}}  + {\displaystyle \frac {1}{2}} \,
{\displaystyle \frac {m - 2}{(m^{2} - 2\,m + 1)\,n}}  \\
 & & \mbox{} - {\displaystyle \frac {3}{2}} \,{\displaystyle 
\frac {1}{(m - 2)\,n}} 
\end{eqnarray*}
\end{maplelatex}

\begin{maplettyout}
covering   rows   1   2   3
\end{maplettyout}

\begin{maplettyout}
using   cases   62   1   3   7
\end{maplettyout}

\begin{maplettyout}
\begin{center}{\bf Case   68}\end{center}
\end{maplettyout}

\begin{maplelatex}
\begin{eqnarray*}
\lefteqn{k=5, \,M= \left[ 
{\begin{array}{cccc}
0 & 0 & 0 & 0 \\
{X_{1}} + {\displaystyle \frac {{X_{2}}}{(m - 2)\,(n - 1)}}  & 0
 & {X_{3}} & {X_{4}} \\ [2ex]
{X_{5}} & 0 & {X_{6}} & {X_{7}} \\
0 & {X_{8}} & {X_{9}} & {X_{10}} \\
0 & {X_{11}} & {X_{12}} & {X_{13}}
\end{array}}
 \right] , } \\
 & & F= \left[ 
{\begin{array}{rrrrr}
0 & 0 & 1 & -2 & 1 \\
0 & 1 & 1 & 1 & 0 \\
0 & -2 & 1 & 0 & 0 \\
0 & 1 & 0 & 0 & 0 \\
0 & 0 & 0 & 0 & 0
\end{array}}
 \right]  + 2\,{\displaystyle \frac {2\,m - 5}{(m - 2)\,(m - 1)\,
(m - 3)\,(2\,m\,n - 2\,m - 5\,n + 4)}} 
\end{eqnarray*}
\end{maplelatex}

\begin{maplettyout}
covering   row   1
\end{maplettyout}

\begin{maplettyout}
covering   columns   1   2
\end{maplettyout}

\begin{maplettyout}
using   cases   59   62   3   5   10
\end{maplettyout}

\begin{maplettyout}
\begin{center}{\bf Case   69}\end{center}
\end{maplettyout}

\begin{maplelatex}
\[
k=5, \,M= \left[ 
{\begin{array}{cccc}
0 & 0 & 0 & 0 \\
{X_{1}} & 0 & {X_{2}} & {X_{3}} \\
0 & {X_{4}} & {X_{5}} & {X_{6}} \\
0 & {X_{7}} & {X_{8}} & {X_{9}}
\end{array}}
 \right] , \,F= \left[ 
{\begin{array}{rrrrr}
0 & 0 & 1 & -2 & 1 \\
-1 & -1 & 1 & 1 & 0 \\
3 & -1 & 1 & 0 & 0 \\
-3 & 2 & 0 & 0 & 0 \\
1 & 0 & 0 & 0 & 0
\end{array}}
 \right] 
\]
\end{maplelatex}

\begin{maplettyout}
covering   rows   1   2
\end{maplettyout}

\begin{maplettyout}
covering   column   1
\end{maplettyout}

\begin{maplettyout}
using   cases   65   66   67   68   7   14
\end{maplettyout}

\begin{maplettyout}
\begin{center}{\bf Case   70}\end{center}
\end{maplettyout}

\begin{maplelatex}
\begin{eqnarray*}
\lefteqn{k=5, \,M= \left[ 
{\begin{array}{ccc}
{\displaystyle \frac {{X_{1}}}{(m - 1)\,(n - 1)}}  & 0 & 0 \\
 [2ex]
{X_{2}} & 0 & {X_{3}} \\
0 & {X_{4}} & {X_{5}} \\
0 & {X_{6}} & {X_{7}}
\end{array}}
 \right] , F= \left[ 
{\begin{array}{rrrrr}
0 & -3 & -1 & 1 & 1 \\
-1 & 0 & 0 & 1 & 0 \\
3 & -1 & 1 & 0 & 0 \\
-3 & 2 & 0 & 0 & 0 \\
1 & 0 & 0 & 0 & 0
\end{array}}
 \right]  + 2\,{\displaystyle \frac {m - 2}{(m^{2} - 2\,m + 1)\,(
n - 1)}} } \\
 & & \mbox{} + 2\,{\displaystyle \frac {2\,m - 3}{m\,(m^{2} - 2\,
m + 1)\,(2\,m\,n - 3\,m - 3\,n + 4)}} \mbox{\hspace{162pt}}
\end{eqnarray*}
\end{maplelatex}

\begin{maplettyout}
covering   columns   1   2   3
\end{maplettyout}

\begin{maplettyout}
using   cases   37   57   69   16
\end{maplettyout}

\begin{maplettyout}
\begin{center}{\bf Case   71}\end{center}
\end{maplettyout}

\begin{maplelatex}
\begin{eqnarray*}
\lefteqn{k=5, \,M= \left[ 
{\begin{array}{cccc}
0 & 0 & {X_{1}} & 0 \\
0 & 0 & {X_{2}} & {X_{3}} \\
{X_{4}} + {\displaystyle \frac {{X_{5}}}{(m - 2)\,(n - 1)}}  & {X
_{6}} & 0 & {X_{7}} \\ [2ex]
0 & {X_{8}} & {X_{9}} & {X_{10}}
\end{array}}
 \right] , F= \left[ 
{\begin{array}{rrrrr}
0 & 0 & 1 & -2 & 1 \\
0 & 1 & -3 & 2 & 0 \\
0 & -2 & 2 & 0 & 0 \\
0 & 1 & 0 & 0 & 0 \\
1 & 0 & 0 & 0 & 0
\end{array}}
 \right]  - {\displaystyle \frac {3}{2}} \,{\displaystyle \frac {
1}{(m - 3)\,n}} } \\
 & & \mbox{} + {\displaystyle \frac {1}{2}} \,{\displaystyle 
\frac {m - 3}{(m^{2} - 4\,m + 4)\,n}}  + {\displaystyle \frac {1
}{2}} \,{\displaystyle \frac {(2\,m - 5)^{2}}{(m^{2} - 4\,m + 4)
\,(m - 3)\,(2\,m\,n - 2\,m - 5\,n + 4)}} \mbox{\hspace{23pt}}
\end{eqnarray*}
\end{maplelatex}

\begin{maplettyout}
covering   rows   1   2   3   4
\end{maplettyout}

\begin{maplettyout}
using   cases   59   1   3   7   16
\end{maplettyout}

\begin{maplettyout}
\begin{center}{\bf Case   72}\end{center}
\end{maplettyout}

\begin{maplelatex}
\[
k=5, \,M= \left[ 
{
}
 \right]  - {\displaystyle \frac {1}{2}} \,{\displaystyle \frac {
m - 4}{(m - 2)^{2}\,(n - 2)}} } \\
 & & \mbox{} + 2\,{\displaystyle \frac {m - 2}{(m^{2} - 2\,m + 1)
\,(n - 1)}}  - 2\,{\displaystyle \frac {(2\,m - 3)^{2}}{(m - 2)^{
2}\,(m^{2} - 2\,m + 1)\,m\,(2\,m\,n - 3\,m - 3\,n + 4)}}  \\
 & & \mbox{} - {\displaystyle \frac {3}{2}} \,{\displaystyle 
\frac {1}{m\,(n - 2)}} 
\end{eqnarray*}
\end{maplelatex}

\begin{maplettyout}
covering   columns   1   2   3
\end{maplettyout}

\begin{maplettyout}
using   cases   36   41   7   103   104
\end{maplettyout}

\begin{maplettyout}
\begin{center}{\bf Case   106}\end{center}
\end{maplettyout}

\begin{maplelatex}
\begin{eqnarray*}
\lefteqn{k=4, \,M= \left[ 
{\begin{array}{ccc}
0 & 0 & 0 \\
{X_{1}} & 0 & {X_{2}} \\
0 & {X_{3}} + {\displaystyle \frac {{X_{4}}}{(m - 1)\,n}}  & {X_{
5}}
\end{array}}
 \right] , F= \left[ 
{\begin{array}{rrrr}
0 & 0 & -1 & 1 \\
0 & -1 & 1 & 0 \\
-1 & 1 & 0 & 0 \\
1 & 0 & 0 & 0
\end{array}}
 \right] } \\
 & & \mbox{} - {\displaystyle \frac {(2\,m - 3)^{2}}{(m - 2)^{2}
\,(m^{2} - 2\,m + 1)\,(2\,m\,n - m - 3\,n + 1)}}  - 
{\displaystyle \frac {1}{(m^{2} - 2\,m + 1)\,n}}  \\
 & & \mbox{} + {\displaystyle \frac {1}{(m - 2)^{2}\,(n - 1)}} 
\end{eqnarray*}
\end{maplelatex}

\begin{maplettyout}
covering   rows   1   2   3
\end{maplettyout}

\begin{maplettyout}
using   cases   1   2   5   101
\end{maplettyout}

\begin{maplettyout}
\begin{center}{\bf Case   107}\end{center}
\end{maplettyout}

\begin{maplelatex}
\begin{eqnarray*}
\lefteqn{k=5, \,M= \left[ 
{\begin{array}{cccc}
0 & 0 & 0 & 0 \\
0 & 0 & {X_{1}} + {\displaystyle \frac {{X_{2}}}{(m - 1)\,(n - 1)
}}  & {X_{3}} \\ [2ex]
{X_{4}} & {X_{5}} & 0 & {X_{6}} \\
0 & {X_{7}} & {X_{8}} & {X_{9}}
\end{array}}
 \right] , F= \left[ 
{\begin{array}{rrrrr}
0 & 0 & 0 & -1 & 1 \\
0 & 0 & -1 & 1 & 0 \\
1 & -2 & 1 & 0 & 0 \\
-2 & 2 & 0 & 0 & 0 \\
1 & 0 & 0 & 0 & 0
\end{array}}
 \right] } \\
 & & \mbox{} + {\displaystyle \frac {1}{(m^{2} - 4\,m + 4)\,(n - 
2)}}  - {\displaystyle \frac {1}{(m - 1)^{2}\,(n - 1)}}  \\
 & & \mbox{} - {\displaystyle \frac {(2\,m - 3)^{2}}{(m - 1)^{2}
\,(m^{2} - 4\,m + 4)\,(2\,m\,n - 3\,m - 3\,n + 4)}} 
\mbox{\hspace{62pt}}
\end{eqnarray*}
\end{maplelatex}

\begin{maplettyout}
covering   rows   1   2   3   4
\end{maplettyout}

\begin{maplettyout}
using   cases   1   2   10   101   106
\end{maplettyout}

\begin{maplettyout}
\begin{center}{\bf Case   108}\end{center}
\end{maplettyout}

\begin{maplelatex}
\begin{eqnarray*}
\lefteqn{k=4, \,M= \left[ 
{\begin{array}{rcc}
0 & 0 & 0 \\
0 & 0 & {X_{1}} \\
0 & {X_{2}} + {\displaystyle \frac {{X_{3}}}{(m - 1)\,n}}  & {X_{
4}}
\end{array}}
 \right] , F= \left[ 
{\begin{array}{rrrr}
0 & 0 & -1 & 1 \\
-1 & 0 & 1 & 0 \\
0 & 0 & 0 & 0 \\
1 & 0 & 0 & 0
\end{array}}
 \right] } \\
 & & \mbox{} - {\displaystyle \frac {(2\,m - 3)^{2}}{(m^{4} - 6\,
m^{3} + 13\,m^{2} - 12\,m + 4)\,(2\,m\,n - m - 3\,n + 1)}}  - 
{\displaystyle \frac {1}{(m^{2} - 2\,m + 1)\,n}}  \\
 & & \mbox{} + {\displaystyle \frac {1}{(m^{2} - 4\,m + 4)\,(n - 
1)}} 
\end{eqnarray*}
\end{maplelatex}

\begin{maplettyout}
covering   rows   1   2   3
\end{maplettyout}

\begin{maplettyout}
using   cases   1   2   3   101
\end{maplettyout}

\begin{maplettyout}
\begin{center}{\bf Case   109}\end{center}
\end{maplettyout}

\begin{maplelatex}
\begin{eqnarray*}
\lefteqn{k=5, \,M= \left[ 
{\begin{array}{rccc}
0 & 0 & 0 & 0 \\
0 & 0 & 0 & {X_{1}} \\
0 & 0 & {X_{2}} + {\displaystyle \frac {{X_{3}}}{(m - 1)\,(n - 1)
}}  & {X_{4}} \\ [2ex]
0 & {X_{5}} & {X_{6}} & {X_{7}}
\end{array}}
 \right] , F= \left[ 
{\begin{array}{rrrrr}
0 & 0 & 0 & -1 & 1 \\
0 & -1 & 0 & 1 & 0 \\
0 & 0 & 0 & 0 & 0 \\
-1 & 1 & 0 & 0 & 0 \\
1 & 0 & 0 & 0 & 0
\end{array}}
 \right] } \\
 & & \mbox{} - {\displaystyle \frac {(2\,m - 3)^{2}}{(m - 1)^{2}
\,(m - 2)^{2}\,(2\,m\,n - 3\,m - 3\,n + 4)}}  - {\displaystyle 
\frac {1}{(m - 1)^{2}\,(n - 1)}}  + {\displaystyle \frac {1}{(m
 - 2)^{2}\,(n - 2)}} 
\end{eqnarray*}
\end{maplelatex}

\begin{maplettyout}
covering   rows   1   2   3   4
\end{maplettyout}

\begin{maplettyout}
using   cases   1   2   3   101   108
\end{maplettyout}

\begin{maplettyout}
\begin{center}{\bf Case   110}\end{center}
\end{maplettyout}

\begin{maplelatex}
\begin{eqnarray*}
\lefteqn{k=5, \,M= \left[ 
{\begin{array}{ccc}
{X_{1}} + {\displaystyle \frac {{X_{2}}}{(m - 1)\,(n - 1)}}  & 0
 & 0 \\ [2ex]
0 & 0 & {X_{3}} \\
{X_{4}} & {X_{5}} & 0 \\
0 & {X_{6}} & {X_{7}}
\end{array}}
 \right] , F= \left[ 
{\begin{array}{rrrrr}
0 & -1 & 0 & 1 & 1 \\
0 & 0 & -4 & 1 & 0 \\
1 & -3 & 2 & 0 & 0 \\
-2 & 2 & 0 & 0 & 0 \\
1 & 0 & 0 & 0 & 0
\end{array}}
 \right]  - {\displaystyle \frac {1}{2}} \,{\displaystyle \frac {
m - 4}{(m - 2)^{2}\,(n - 2)}} } \\
 & & \mbox{} + 2\,{\displaystyle \frac {m - 2}{(m^{2} - 2\,m + 1)
\,(n - 1)}}  - 2\,{\displaystyle \frac {(2\,m - 3)^{2}}{(m^{2} - 
2\,m + 1)\,m\,(m - 2)^{2}\,(2\,m\,n - 3\,m - 3\,n + 4)}}  \\
 & & \mbox{} + {\displaystyle \frac {5}{2}} \,{\displaystyle 
\frac {1}{m\,(n - 2)}} 
\end{eqnarray*}
\end{maplelatex}

\begin{maplettyout}
covering   columns   1   2   3
\end{maplettyout}

\begin{maplettyout}
using   cases   47   3   26   107   109
\end{maplettyout}

\begin{maplettyout}
\begin{center}{\bf Case   111}\end{center}
\end{maplettyout}

\begin{maplelatex}
\begin{eqnarray*}
\lefteqn{k=5, \,M= \left[ 
{\begin{array}{ccc}
{X_{1}} + {\displaystyle \frac {{X_{2}}}{(m - 1)\,(n - 1)}}  & 0
 & 0 \\ [2ex]
0 & 0 & {X_{3}} \\
0 & {X_{4}} & {X_{5}}
\end{array}}
 \right] , F= \left[ 
{\begin{array}{rrrrr}
0 & -4 & 0 & 1 & 1 \\
-4 & 3 & 0 & 1 & 0 \\
0 & 0 & 0 & 0 & 0 \\
1 & 1 & 0 & 0 & 0 \\
1 & 0 & 0 & 0 & 0
\end{array}}
 \right]  + {\displaystyle \frac {1}{(m - 2)^{2}\,(n - 2)}} } \\
 & & \mbox{} + {\displaystyle \frac {27}{4}} \,{\displaystyle 
\frac {1}{(3\,m - 4)\,n}}  - 4\,{\displaystyle \frac {1}{(m^{2}
 - 2\,m + 1)\,(n - 1)}}  \\
 & & \mbox{} - 4\,{\displaystyle \frac {(4\,m^{2} - 12\,m + 9)\,(
2\,m - 3)}{(m^{2} - 2\,m + 1)\,(m - 2)^{2}\,m\,(3\,m - 4)\,(2\,m
\,n - 3\,m - 3\,n + 4)}}  \\
 & & \mbox{} + {\displaystyle \frac {1}{2}} \,{\displaystyle 
\frac {1}{(m - 2)\,(n - 2)}}  - {\displaystyle \frac {1}{2}} \,
{\displaystyle \frac {1}{(m - 2)\,n}}  + {\displaystyle \frac {1
}{4}} \,{\displaystyle \frac {1}{m\,n}}  - {\displaystyle \frac {
1}{2}} \,{\displaystyle \frac {1}{m\,(n - 2)}} 
\mbox{\hspace{87pt}}
\end{eqnarray*}
\end{maplelatex}

\begin{maplettyout}
covering   rows   1   2   3
\end{maplettyout}

\begin{maplettyout}
using   cases   39   53   105   110
\end{maplettyout}

\begin{maplettyout}
\begin{center}{\bf Case   112}\end{center}
\end{maplettyout}

\begin{maplelatex}
\begin{eqnarray*}
\lefteqn{k=5, \,M= \left[ 
{\begin{array}{cccc}
{X_{1}} + {\displaystyle \frac {{X_{2}}}{(m - 1)\,(n - 1)}}  & 0
 & {X_{3}} & 0 \\ [2ex]
0 & {X_{4}} & 0 & {X_{5}} \\
{X_{6}} & 0 & {X_{7}} & {X_{8}} \\
0 & {X_{9}} & {X_{10}} & {X_{11}}
\end{array}}
 \right] , F= \left[ 
{\begin{array}{rrrrr}
0 & 0 & 1 & -2 & 1 \\
0 & 1 & -3 & 2 & 0 \\
1 & -3 & 2 & 0 & 0 \\
-2 & 2 & 0 & 0 & 0 \\
1 & 0 & 0 & 0 & 0
\end{array}}
 \right] } \\
 & & \mbox{} + {\displaystyle \frac {1}{(m^{2} - 4\,m + 4)\,(n - 
2)}}  - {\displaystyle \frac {1}{(m^{2} - 2\,m + 1)\,(n - 1)}} 
 \\
 & & \mbox{} - {\displaystyle \frac {(2\,m - 3)^{2}}{(m^{2} - 2\,
m + 1)\,(m^{2} - 4\,m + 4)\,(2\,m\,n - 3\,m - 3\,n + 4)}} 
\mbox{\hspace{40pt}}
\end{eqnarray*}
\end{maplelatex}

\begin{maplettyout}
covering   rows   1   2   3   4
\end{maplettyout}

\begin{maplettyout}
using   cases   1   5   16   101   102
\end{maplettyout}

\begin{maplettyout}
\begin{center}{\bf Case   113}\end{center}
\end{maplettyout}

\begin{maplelatex}
\begin{eqnarray*}
\lefteqn{k=5, \,M= \left[ 
{\begin{array}{cccc}
{X_{1}} + {\displaystyle \frac {{X_{2}}}{(m - 1)\,(n - 1)}}  & {X
_{3}} & 0 & 0 \\ [2ex]
0 & 0 & {X_{4}} & {X_{5}} \\
0 & {X_{6}} & {X_{7}} & {X_{8}}
\end{array}}
 \right] , F= \left[ 
{\begin{array}{rrrrr}
0 & -1 & 3 & -3 & 1 \\
-3 & 2 & -1 & 2 & 0 \\
0 & 0 & 0 & 0 & 0 \\
1 & 1 & 0 & 0 & 0 \\
1 & 0 & 0 & 0 & 0
\end{array}}
 \right] } \\
 & & \mbox{} + {\displaystyle \frac {1}{2}} \,{\displaystyle 
\frac {1}{(m - 2)\,(n - 2)}}  - {\displaystyle \frac {1}{2}} \,
{\displaystyle \frac {1}{(m - 2)\,n}}  + {\displaystyle \frac {1
}{2}} \,{\displaystyle \frac {m}{(m - 2)^{2}\,(n - 2)}}  - 2\,
{\displaystyle \frac {m}{(m - 1)^{2}\,(n - 1)}}  \\
 & & \mbox{} + {\displaystyle \frac {9}{2}} \,{\displaystyle 
\frac {1}{(3\,m - 4)\,n}}  - 2\,{\displaystyle \frac {(4\,m^{2}
 - 12\,m + 9)\,(2\,m - 3)}{(3\,m - 4)\,(m - 1)^{2}\,(m - 2)^{2}\,
(2\,m\,n - 3\,m - 3\,n + 4)}} 
\end{eqnarray*}
\end{maplelatex}

\begin{maplettyout}
covering   rows   1   2   3
\end{maplettyout}

\begin{maplettyout}
using   cases   41   56   16   112
\end{maplettyout}

\begin{maplettyout}
\begin{center}{\bf Case   114}\end{center}
\end{maplettyout}

\begin{maplelatex}
\begin{eqnarray*}
\lefteqn{k=5, \,M= \left[ 
{\begin{array}{ccc}
{X_{1}} + {\displaystyle \frac {{X_{2}}}{(m - 1)\,(n - 1)}}  & 0
 & 0 \\ [2ex]
{X_{3}} & 0 & {X_{4}} \\
0 & {X_{5}} & {X_{6}} \\
0 & {X_{7}} & {X_{8}}
\end{array}}
 \right] , F= \left[ 
{\begin{array}{rrrrr}
0 & -3 & 0 & 1 & 1 \\
-1 & 0 & 0 & 1 & 0 \\
3 & -1 & 1 & 0 & 0 \\
-3 & 2 & 0 & 0 & 0 \\
1 & 0 & 0 & 0 & 0
\end{array}}
 \right] } \\
 & & \mbox{} - {\displaystyle \frac {1}{2}} \,{\displaystyle 
\frac {m - 4}{(m^{2} - 4\,m + 4)\,(n - 2)}}  + 2\,{\displaystyle 
\frac {m - 2}{(m - 1)^{2}\,(n - 1)}}  \\
 & & \mbox{} - 2\,{\displaystyle \frac {(2\,m - 3)^{2}}{(m^{2} - 
4\,m + 4)\,(m - 1)^{2}\,m\,(2\,m\,n - 3\,m - 3\,n + 4)}}  - 
{\displaystyle \frac {1}{2}} \,{\displaystyle \frac {1}{m\,(n - 2
)}} 
\end{eqnarray*}
\end{maplelatex}

\begin{maplettyout}
covering   columns   1   2   3
\end{maplettyout}

\begin{maplettyout}
using   cases   37   69   16   112
\end{maplettyout}

\begin{maplettyout}
\begin{center}{\bf Case   115}\end{center}
\end{maplettyout}

\begin{maplelatex}
\[
k=5, \,M= \left[ 
{\begin{array}{ccc}
{X_{1}} & 0 & 0 \\
0 & {X_{2}} & {X_{3}} \\
0 & {X_{4}} & {X_{5}}
\end{array}}
 \right] , \,F= \left[ 
{\begin{array}{rrrrr}
1 & 1 & -4 & 1 & 1 \\
1 & 5 & -1 & 1 & 0 \\
-4 & -1 & 2 & 0 & 0 \\
1 & 1 & 0 & 0 & 0 \\
1 & 0 & 0 & 0 & 0
\end{array}}
 \right] 
\]
\end{maplelatex}

\begin{maplettyout}
covering   row   1
\end{maplettyout}

\begin{maplettyout}
covering   column   1
\end{maplettyout}

\begin{maplettyout}
using   cases   72   111   113   114
\end{maplettyout}

\begin{maplettyout}
\begin{center}{\bf Case   116}\end{center}
\end{maplettyout}

\begin{maplelatex}
\[
k=5, \,M= \left[ 
{\begin{array}{ccc}
{X_{1}} & {X_{2}} & 0 \\
{X_{3}} & 0 & {X_{4}} \\
0 & {X_{5}} & {X_{6}} \\
0 & {X_{7}} & {X_{8}}
\end{array}}
 \right] , \,F= \left[ 
{\begin{array}{rrrrr}
1 & 1 & -4 & 1 & 1 \\
-4 & 1 & 2 & 1 & 0 \\
6 & -5 & 2 & 0 & 0 \\
-4 & 3 & 0 & 0 & 0 \\
1 & 0 & 0 & 0 & 0
\end{array}}
 \right] 
\]
\end{maplelatex}

\begin{maplettyout}
covering   columns   1   2   3
\end{maplettyout}

\begin{maplettyout}
using   cases   72   20   98
\end{maplettyout}

\begin{maplettyout}
\begin{center}{\bf Case   117}\end{center}
\end{maplettyout}

\begin{maplelatex}
\[
k=5, \,M= \left[ 
{\begin{array}{cc}
{X_{1}} & 0 \\
0 & {X_{2}} \\
0 & {X_{3}}
\end{array}}
 \right] , \,F= \left[ 
{\begin{array}{rrrrr}
-4 & 1 & 1 & 1 & 1 \\
11 & -7 & 1 & 1 & 0 \\
-9 & 5 & 1 & 0 & 0 \\
1 & 1 & 0 & 0 & 0 \\
1 & 0 & 0 & 0 & 0
\end{array}}
 \right] 
\]
\end{maplelatex}

\begin{maplettyout}
covering   columns   1   2
\end{maplettyout}

\begin{maplettyout}
using   cases   99   115   116
\end{maplettyout}

\begin{maplettyout}
\begin{center}{\bf Case   118}\end{center}
\end{maplettyout}

\begin{maplelatex}
\[
k=5, \,M= \left[ 
{\begin{array}{r}
0 \\
0
\end{array}}
 \right] , \,F= \left[ 
{\begin{array}{rrrrr}
6 & 1 & 1 & 1 & 1 \\
-9 & 1 & 1 & 1 & 0 \\
1 & 1 & 1 & 0 & 0 \\
1 & 1 & 0 & 0 & 0 \\
1 & 0 & 0 & 0 & 0
\end{array}}
 \right] 
\]
\end{maplelatex}

\begin{maplettyout}
covering   column   1
\end{maplettyout}

\begin{maplettyout}
using   cases   100   117
\end{maplettyout}

\begin{maplettyout}
\begin{center}{\bf Case   119}\end{center}
\end{maplettyout}

\begin{maplelatex}
\[
k=5, \,M= \left[ 
{\begin{array}{cccc}
{X_{1}} & {X_{2}} & {X_{3}} & 0 \\
{X_{4}} & {X_{5}} & 0 & {X_{6}} \\
{X_{7}} & 0 & {X_{8}} & {X_{9}} \\
0 & {X_{10}} & {X_{11}} & {X_{12}}
\end{array}}
 \right] , \,F= \left[ 
{\begin{array}{rrrrr}
1 & -4 & 6 & -4 & 1 \\
-4 & 12 & -12 & 4 & 0 \\
6 & -12 & 6 & 0 & 0 \\
-4 & 4 & 0 & 0 & 0 \\
1 & 0 & 0 & 0 & 0
\end{array}}
 \right] 
\]
\end{maplelatex}

\begin{maplettyout}
covering   rows   1   2   3   4
\end{maplettyout}

\begin{maplettyout}
using   cases   1   20
\end{maplettyout}

\begin{maplettyout}
\begin{center}{\bf Case   120}\end{center}
\end{maplettyout}

\begin{maplelatex}
\[
k=5, \,M= \left[ 
{\begin{array}{ccc}
{X_{1}} & {X_{2}} & 0 \\
{X_{3}} & 0 & {X_{4}} \\
0 & {X_{5}} & {X_{6}}
\end{array}}
 \right] , \,F= \left[ 
{\begin{array}{rrrrr}
-4 & 11 & -9 & 1 & 1 \\
11 & -21 & 9 & 1 & 0 \\
-9 & 9 & 0 & 0 & 0 \\
1 & 1 & 0 & 0 & 0 \\
1 & 0 & 0 & 0 & 0
\end{array}}
 \right] 
\]
\end{maplelatex}

\begin{maplettyout}
covering   rows   1   2   3
\end{maplettyout}

\begin{maplettyout}
using   cases   116   119
\end{maplettyout}

\begin{maplettyout}
\begin{center}{\bf Case   121}\end{center}
\end{maplettyout}

\begin{maplelatex}
\[
k=5, \,M= \left[ 
{\begin{array}{cc}
{X_{1}} & 0 \\
0 & {X_{2}}
\end{array}}
 \right] , \,F= \left[ 
{\begin{array}{rrrrr}
6 & -9 & 1 & 1 & 1 \\
-9 & 7 & 1 & 1 & 0 \\
1 & 1 & 1 & 0 & 0 \\
1 & 1 & 0 & 0 & 0 \\
1 & 0 & 0 & 0 & 0
\end{array}}
 \right] 
\]
\end{maplelatex}

\begin{maplettyout}
covering   rows   1   2
\end{maplettyout}

\begin{maplettyout}
using   cases   117   120
\end{maplettyout}

\begin{maplettyout}
\begin{center}{\bf Case   122}\end{center}
\end{maplettyout}

\begin{maplelatex}
\[
k=5, \,M= \left[ 
{\begin{array}{r}
0
\end{array}}
 \right] , \,F= \left[ 
{\begin{array}{rrrrr}
-4 & 1 & 1 & 1 & 1 \\
1 & 1 & 1 & 1 & 0 \\
1 & 1 & 1 & 0 & 0 \\
1 & 1 & 0 & 0 & 0 \\
1 & 0 & 0 & 0 & 0
\end{array}}
 \right] 
\]
\end{maplelatex}

\begin{maplettyout}
covering   row   1
\end{maplettyout}

\begin{maplettyout}
using   cases   118   121
\end{maplettyout}

\begin{maplettyout}
\begin{center}{\bf Case   123}\end{center}
\end{maplettyout}

\begin{maplelatex}
\[
k=5, \,M=[], \,F= \left[ 
{\begin{array}{rrrrr}
1 & 1 & 1 & 1 & 1 \\
1 & 1 & 1 & 1 & 0 \\
1 & 1 & 1 & 0 & 0 \\
1 & 1 & 0 & 0 & 0 \\
1 & 0 & 0 & 0 & 0
\end{array}}
 \right] 
\]
\end{maplelatex}

\begin{maplettyout}
using   case   122
\end{maplettyout}
